\documentclass[english]{article}
\usepackage[T1]{fontenc}
\usepackage[latin9]{inputenc}
\usepackage[letterpaper]{geometry}
\usepackage{float}
\usepackage{amsmath}
\usepackage{graphicx}
\usepackage{color}

\newtheorem{theorem}{Theorem}
\newtheorem{assumption}{Assumption}
\newtheorem{lemma}{Lemma}

\usepackage{amsfonts} 
\makeatother

\usepackage{babel}

\title{Generalized multiscale finite element methods for wave propagation
in heterogeneous media}
\author{Eric T. Chung\thanks{Department of Mathematics, The Chinese University of Hong Kong, Hong Kong SAR. Eric Chung's research
is supported by the Hong Kong RGC General Research Fund project 400411.},
Yalchin Efendiev\thanks{Department of Mathematics, Texas A\&M University, College Station, TX 77843, USA \& center for Numerical Porous Media (NumPor),
King Abdullah University of Science and Technology (KAUST),
Thuwal 23955-6900, Kingdom of Saudi Arabia}\,
 and Wing Tat Leung\thanks{Department of Mathematics, The Chinese University of Hong Kong, Hong Kong SAR.}
}

\begin{document}

\maketitle

\begin{abstract}
Numerical modeling of wave propagation in heterogeneous media is important
in many applications. Due to the complex nature, direct
numerical simulations on the fine grid are prohibitively expensive. It is therefore
important to develop efficient and accurate methods that allow the
use of coarse grids. In this paper, we present a multiscale finite
element method for  wave propagation on a coarse grid.
The proposed method is based on the Generalized Multiscale Finite Element
Method (GMsFEM) (see \cite{egh12}).
To construct
multiscale basis functions, we start with two snapshot spaces in each
coarse-grid block where
one represents the degrees of freedom on the boundary and the other
represents the degrees of freedom in the interior.
 We use local spectral problems to identify important
modes in each snapshot space.
These local spectral problems are different from each other and their
formulations are based on the analysis. To our best knowledge,
this is the first time where
multiple snapshot spaces
and multiple spectral problems are used and
necessary for efficient computations.
Using the dominant modes from local spectral problems,
multiscale basis functions are constructed to
represent the solution space locally within each coarse block.
These multiscale basis functions
are coupled via the symmetric interior penalty discontinuous Galerkin
method which provides a block diagonal mass matrix, and, consequently,
results in fast computations in an explicit time discretization.
Our methods' stability
and spectral convergence  are rigorously analyzed. Numerical examples
are presented to show our methods' performance. We also test oversampling
strategies. In particular,
we discuss how the modes from different snapshot spaces can affect the proposed methods' accuracy.
\end{abstract}

\section{Introduction}

Numerical modeling of wave propagation is important in many applications
that include geophysics, material science, and so on. For example,
in geophysics applications,  wave propagation simulations play
an important role in determining subsurface properties \cite{virieux:889,virieux:1933,saenger:SM293,GJI:GJI3620,symes:2602,masson:N33}.
These approaches include finite difference methods, finite element
methods, and spectral methods that use polynomials basis \cite{basabe:562,kaser:76,Pelties:2010kx,Hermann:2011fk,GJI:GJI5221,springerlink:10.1007/s00450-010-0109-1,GJI:GJI1653,GJI:GJI967,GJI:GJI4985,CMS,JCAM-lod,nmtma}.
While these methods have different strengths and weaknesses, all of
them tend to have limitations associated with discretization, especially
in 3-D applications as frequency content of the simulated wavefield
increases. Though the solutions to the wave equation have been shown
to be accurate when the grid is fine enough \cite{delprat-jannaud:T37},
the practical limitations in discretization caused by limitations
in computational power restrict this accuracy. An example of an application
where this may be important is in the modeling of fractured media,
where establishing reliable and accurate relationships between the
properties of reflected seismic wavefields and variations in the density,
orientation and compliance of fractures may help provide important
constraints for hydrocarbon production. While more general finite
element and spectral element methods may be able to address some problems
by adapting grids to conform to heterogeneous structures, there are basic
limitations associated with representing fine-scale features,
and there is therefore a need to find approaches that reliably
and accurately incorporate fine-scale features in a coarsely gridded
model.

In this paper, we present a multiscale finite element method for wave
propagation simulations on a coarse grid. The proposed method is based
on the Generalized Multiscale Finite Element Method (GMsFEM) which was
proposed in \cite{egh12} and couples multiscale basis functions via
a discontinuous Galerkin coupling (cf. \cite{eglms13}). To construct
multiscale basis functions, we start with two snapshot spaces in each
coarse-grid block where
one represents the degrees of freedom on the coarse grid's boundary
and the other
represents the degrees of freedom in the interior.
 We use local spectral problems to identify important
modes in each snapshot space.
These local spectral problems are different from each other and their
formulations are based on the analysis.
  Once local basis
functions are identified, we couple these basis functions via Interior
Penalty Discontinuous Galerkin method \cite{IPDG,IPDGbook}.
%This apporach provides a block diagonal
%mass matrix.

%The proposed multiscale finite element methods offer a powerful and
%flexible approach to full waveform modeling that represents the effects
%of fine-scale features on a coarse grid by constructing multiscale
%basis functions. These basis functions, which are computed one time
%using the local fine grid, are then utilized on the coarse grid for
%rapid calculations. The general approach has several potential advantages
%over alternative solutions that instead upscale properties using effective
%medium theories.

Because these basis functions are discontinuous, the interior penalty
discontinuous Galerkin (IPDG) method, for example \cite{IPDG,CiCP,IPDGbook},
is an appropriate choice for solving the time-dependent partial differential
equation. It generally yields a block diagonal mass matrix, hence
the time stepping is very efficient.
The staggered discontinuous Galerkin
methods \cite{ChungEngquist06,ChungEngquist09} have been recently
developed for the accurate wave simulations. By using a carefully
chosen staggered grid, the resulting method is also energy conserving.
Moreover, it is proved that (see \cite{dispersion,JCP-max}) such method
gives smaller dispersion errors, and therefore it is superior for
the wave propagation. The staggered idea has also been extended to other
problems, see for example \cite{OLS,SDG-cd,SDG-curl,JCAM-meta}.
Recently, we have used standard MsFEM basis within staggered methods
\cite{AADA,Geophysics}.
These methods allow some limited upscaling and
provide energy conserving numerical methods on staggered grids.
In this paper, our goal is to construct a systematic
enrichment by appropriately choosing snapshot spaces and corresponding
local spectral problems.

We will focus our discussions on two-dimensional problems. The extension
to three-dimensional
problems is straightforward. Let $\Omega\subset\mathbb{R}^{2}$
be a bounded domain of two dimensions. The paper's aim is to
develop a new multiscale method for the following wave equation \begin{equation}
\frac{\partial^{2}u}{\partial t^{2}}=\nabla\cdot(a\nabla u)+f\quad\quad\mbox{ in }\quad[0,T]\times\Omega\label{eq:waveeqn}\end{equation}
 with the homogeneous Dirichlet boundary condition $u=0$ on $[0,T]\times\partial\Omega$.
The extension to other boundary conditions will be reported in a forthcoming
paper. The function $f(x,t)$ is a given source. The problem (\ref{eq:waveeqn})
is supplemented with the following initial conditions \[
u(x,0)=g_{0}(x),\quad\quad u_{t}(x,0)=g_{1}(x).\]
 We assume that the coefficient $a(x)$ is highly oscillatory, representing
the complicated model in which the waves are simulated. It is well-known
that solving (\ref{eq:waveeqn}) by standard methods requires a very
fine mesh, which is computationally prohibited. Thus a coarse grid
solution strategy is needed. Next we present our fine scale solver.
The fine scale solution is considered as the exact
solution when we discuss the convergence of our multiscale method
 in the following sections. We assume that the domain $\Omega$
is partitioned by a set of rectangles, called fine mesh, with maximum
side length $h>0$. We denote the resulting mesh by $\mathcal{T}^{h}$
and the set of all edges and vertices  by $\mathcal{E}^{h}$
and $\mathcal{N}^{h}$ respectively. We assume that the fine-mesh
discretization of the wave equation provides an accurate approximation
of the solution. The fine scale solver is the standard conforming
bilinear finite element method. Let $V_{h}$ be the standard conforming piecewise bilinear finite
element space. We find $u_{h}\in V_{h}$ such that \begin{equation}
(\frac{\partial^{2}u_{h}}{\partial t^{2}},v)+a(u_{h},v)=(f,v),\quad\forall v\in V_{h},\label{eq:globalFEM}\end{equation}
 where the bilinear form $a$ is defined by \begin{equation}
a(u,v)=\int_{\Omega}a\nabla u\cdot\nabla v,\quad\forall u,v\in V_{h}\end{equation}
 and $(\cdot,\cdot)$ represents the standard $L^{2}$ inner product
defined on $\Omega$.

The numerical results are presented for several representative
examples. We investigate the GMsFEM's accuracy  and, in particular,
how choosing modes from different snapshot spaces can affect the
accuracy. Our numerical results show that choosing the basis functions
from interior modes can improve the accuracy of GMsFEM substantially
for wave equations. These results differ from those we observe for
flow equations \cite{egh12}.

The paper is organized as follows. In Section \ref{sec:method}, we
will present the new multiscale method. Numerical results are shown
in Section \ref{sec:num}. Stability and spectral convergence of the semi-discrete
scheme are proved in Section \ref{sec:stab}. In Section \ref{sec:full},
the convergence of the fully-discrete scheme is also proved. Finally,
conclusions are presented.

\section{The generalized multiscale finite element method}

\label{sec:method}

In this section, we will give a detailed description of our new generalized
multiscale finite element method. The method gives a numerical solver
on a coarse grid, providing an efficient way to simulate waves in complicated
media. As we will discuss next, the local basis functions are obtained
via the solutions of some local spectral problems which are used to
obtain the most dominant modes. These modes form the basis functions
of our multiscale finite element method.

We introduce a coarse mesh that consists of union of connected fine-mesh
grid blocks which is denoted by $\mathcal{T}^{H}$ and the set of all
edges  by $\mathcal{E}^{H}$. We denote the size of the
coarse mesh by $H$. Even though it is convenient to choose rectangular
coarse grid blocks, the shapes of the coarse grid blocks can be quite
general and our analysis can be applied without the assumption of
rectangular coarse grid blocks.

For each coarse grid block $K$, we define $\partial\mathcal{T}^{h}(K)$
be the restriction of the conforming piecewise bilinear functions with respect
to the fine mesh on $\partial K$.
We remark that, for a coarse grid edge $e\in\mathcal{E}^H$ that is shared by two coarse grid blocks $K_1$ and $K_2$,
the values of the two functions in $\partial\mathcal{T}^{h}(K_1)$ and $\partial\mathcal{T}^{h}(K_2)$
on $e$ are in general different.
The union of all $\partial\mathcal{T}^{h}(K)$
is denoted by $\partial\mathcal{T}^{h}$. Moreover, we define $H^{1}(\mathcal{T}^{H})$
as the space of functions whose restrictions on $K$ belongs to $H^{1}(K)$.

\subsection{Global IPDG solver}

We will apply the standard symmetric IPDG approach to solve (\ref{eq:waveeqn})
on the coarse grid $\mathcal{T}^H$. The method follows the standard framework as discussed
in \cite{IPDG,IPDGbook}, but the finite element space will be replaced
by the space spanned by our multiscale basis functions. We emphasize
that the use of the
IPDG approach is an example of the global coupling of our local multiscale
basis functions, and other choices of coarse grid methods are equally
good. The key to  our proposed method's success of is the construction
of our local multiscale basis functions.

First, we introduce some notations. For each interior coarse edge
$e\in\mathcal{E}^{H}$, we let $K^{-}$ and $K^{+}$ be the two coarse
grid blocks having the common coarse edge $e$. Then we define the
average and the jump operators respectively by \begin{align*}
\{v\}_{e} & =\frac{v^{+}+v^{-}}{2},\\
{}[u]_{e} & =u^{+}-u^{-},\end{align*}
 where $u^{\pm}=u|_{K^{\pm}}$ and we have assumed that the normal
vector on $e$ is pointing from $K^{+}$ to $K^{-}$. For each coarse
edge $e$ that lies on the boundary of $\Omega$, we define \[
\{v\}_{e}=v,\quad[u]_{e}=u\]
 assuming the unit normal vector on $e$ is pointing outside the domain.
Let $V_{H}$ be a finite dimensional function space which consists
of functions that are smooth on each coarse grid blocks but are in
general discontinuous across coarse grid edges. We can then state
the IPDG method as: find $u_{H}(t,\cdot)\in V_{H}$ such that \begin{equation}
(\frac{\partial u_{H}}{\partial t^{2}},v)+a_{DG}(u_{H},v)=l(v),\quad\quad\forall\, v\in V_{H},\label{eq:ipdg}\end{equation}
 where the bilinear form $a_{DG}(u,v)$ and the linear functional
$l(v)$ are defined by \begin{align*}
a_{DG}(u,v) & =\sum_{K\in\mathcal{T}^{H}}\int_{K}a\nabla u\cdot\nabla v+\sum_{e\in\mathcal{E}^{H}}\Big(-\int_{e}\{a\nabla u\cdot n\}_{e}\,[v]_{e}-\int_{e}\{a\nabla v\cdot n\}_{e}\,[u]_{e}+\cfrac{\gamma}{h}\int_{e}a[u]_{e}\,[v]_{e}\Big)\\
l(v) & =(f,v)\end{align*}
 where $\gamma>0$ is a penalty parameter and $n$ denotes the unit
normal vector on $e$. The initial conditions for the problem (\ref{eq:ipdg})
are defined by $u_{H}(0)=P_{H}(g_{0})$ and $(u_{H})_{t}(0)=P_{H}(g_{1})$,
where $P_{H}$ is the $L^{2}$-projection operator into $V_{H}$.

Let $T>0$ be a fixed time and $\Delta t=T/N$ be the time step size.
The time discretization is done in the standard way, we find $u_{H}^{n+1}\in V_{H}$
such that \begin{equation}
(u_{H}^{n+1},v)=2(u_{H}^{n},v)-(u_{H}^{n-1},v)-\Delta t^{2}\Big(a_{DG}(u_{H}^{n},v)-l(v)\Big),\quad\quad\forall\, v\in V_{H}\label{eq:fulldiscrete}\end{equation}
 in each time step. Throughout the paper, the notation $u^{n}$ represents
the value of the function $u$ at time $t_{n}$. The initial conditions
are obtained as follows \begin{align*}
u_{H}^{0} & =P_{H}(g_{0}),\\
u_{H}^{1} & =u_{H}^{0}+\Delta t\, P_{H}(g_{1})+\cfrac{\Delta t^{2}}{2}\,\tilde{v},\end{align*}
 where $\tilde{v}\in V_{H}$ is defined by \[
(\tilde{v},v)=(f(0),v)-a_{DG}(g_{0},v),\quad\;\forall\, v\in V_{H}.\]

\subsection{Multiscale basis functions}

In this section, we will give the definition of the space $V_{H}$.
We will discuss the choice
of our basis functions on one single
coarse grid block $K$, as the definitions on other coarse grid blocks are similar.
We recall that $K$ is the union of a set of rectangular elements.
We decompose the space $V_{H}$ into two components,
namely
\[
V_{H}=V_{H}^{1}+V_{H}^{2}.
\]
The restrictions of $V_{H},V_{H}^{1}$ and $V_{H}^{2}$ on $K$ are
denoted by $V_{H}(K)$, $V_{H}^{1}(K)$, and $V_{H}^{2}(K)$ respectively.
Moreover, the restriction of the conforming space $V_{h}$ on $K$
is denoted by $V_{h}(K)$.

{\bf Definition of $V_{H}^{1}(K)$.} To define $V_{H}^{1}(K)$, for each fine
grid node $x_{i}$ on the boundary of $K$, we find $w_{i,K}\in V_{h}(K)$
by solving \begin{equation}
\int_{K}a\nabla w_{i,K}\cdot\nabla v=0,\quad\quad\forall\, v\in V_{h}(K)\label{eq:homo extention}\end{equation}
 with boundary condition $w_{i,K}=1$ at $x_{i}$ and $w_{i,K}=0$
at the other grid points on the boundary of $K$. The functions
$w_{i,K}$ defined above are the $a$-harmonic extensions of the unit
basis functions of $\partial\mathcal{T}^{h}(K)$.
We let $n$ be the number of these $a$-harmonic extensions
and define \[
V_{H}^{1}(K)=\text{span}\{w_{1,K},\cdots,w_{n,K}\}.\]
We remark that $n$ is the number of boundary grid points on $\partial K$
and its value changes with $K$.
 In our numerical simulations, we do not need to use all of these
basis functions and use a local spectral problem in the space of snapshots
to identify multiscale basis functions by choosing dominant modes.
We use $E$ to denote the sum of the reciprocal of
the eigenvalues of the local spectral
problem. We will choose the eigenfunctions corresponding to small
eigenvalues so that the sum of the reciprocals
these small eigenvalues is a small
percentage of $E$.
The use of eigenfunctions corresponding to small eigenvalues
means that we use the coarse component in $V_H^1(K)$
as the approximation space.

{\bf Local spectral problem on $V_{H}^{1}(K)$.} The spectral problem we propose is \begin{equation}
\int_{K}a\nabla w_{\mu}\cdot\nabla v=\cfrac{\mu}{H}\int_{\partial K}w_{\mu}v,\quad\quad\forall\, v\in V_{H}^{1}(K).\label{eq:eigenspace1}\end{equation}
 We assume that the eigenvalues are ordered so that $0=\mu_{1}< \mu_{2}\leq \mu_3 \leq \cdots\leq\mu_{n}$.
The corresponding eigenfunctions in the snapshot space are denoted by $\widetilde{w}_{i,K}$,
$i=1,2,\cdots,n$, which are normalized with respect to the $L^2$-norm on $\partial K$.
We denote the total energy on the coarse grid
block $K$ by $E_{K}$ which is defined by $E_K = \sum_{i=2}^n \mu_i^{-1}$.
We can then choose the first $p$ eigenvalues
so that the sum $\sum_{i=2}^{p} \mu_{i}^{-1}$ %
is a portion of the total energy $E_{K}$.
Note that we can take different $p$ for different coarse grid blocks $K$.
Finally we define \[
\widetilde{V}_{H}^{1}(K)=\text{span}\{\widetilde{w}_{1,K},\cdots,\widetilde{w}_{p,K}\}.\]
 Clearly, $\widetilde{V}_{H}^{1}(K)\subset V_{H}^{1}(K)$.

{\bf Definition of $V_{H}^{2}(K)$.} The space $V_{H}^{2}(K)$ contains functions in $V_h(K)$ that are zero on the
boundary of $K$, which is denoted by $V_h^0(K)$.

{\bf Local spectral problem on $V_{H}^{2}(K)$.}
We will use a suitable (another) spectral problem  to identify
the important modes. The proposed eigenvalue problem has the
following form: find $z_{\lambda}\in V_{h}^0(K)$
such that \begin{equation}
\int_{K}a\nabla z_{\lambda}\cdot\nabla v=\cfrac{\lambda}{H^{2}}\int_{K}z_{\lambda}v,\quad\quad\forall\, v\in V_{h}^0(K).\label{eq:eigenproblem}\end{equation}
 Assume that the eigenvalues are ordered so that $\lambda_{1,K}\leq\lambda_{2,K}\leq\cdots$
and the corresponding eigenfunctions are denoted by $z_{i,K}$, which are normalized with respect to the $L^2$-norm on $K$.
In
practice, for each coarse grid block, we can take the first $m$ eigenfunctions,
and the space $V_{H}^{2}(K)$ is spanned by these functions, that
is \[
V_{H}^{2}(K)=\text{span}\{z_{1,K},\cdots,z_{m,K}\}.\]
 In principle, one can choose different numbers of eigenfunctions
for the space $V_{H}^{2}(K)$ for different coarse grid blocks. Nevertheless,
our numerical results show that only the first few eigenfunctions
are enough to obtain a reasonable accuracy.
% and there is no need to
%choose different $m$ for different coarse grid blocks.

{\bf Orthogonality of $V_{H}^{1}(K)$ and  $V_{H}^{2}(K)$.}  Finally, we point out the following orthogonality condition which
will be used in our analysis. For any $v\in V_{H}^{1}(K)$ and $u\in V_{H}^{2}(K)$,
we conclude by (\ref{eq:homo extention}) that \begin{equation}
\int_{K}a\nabla v\cdot\nabla u=0.\label{eq:Orthogonal}\end{equation}
 This means that the two spaces $V_{H}^{1}(K)$ and $V_{H}^{2}(K)$
are orthogonal.

\section{Numerical Results}

\label{sec:num}

\begin{figure}[H]
 \centering \includegraphics[scale=0.5]{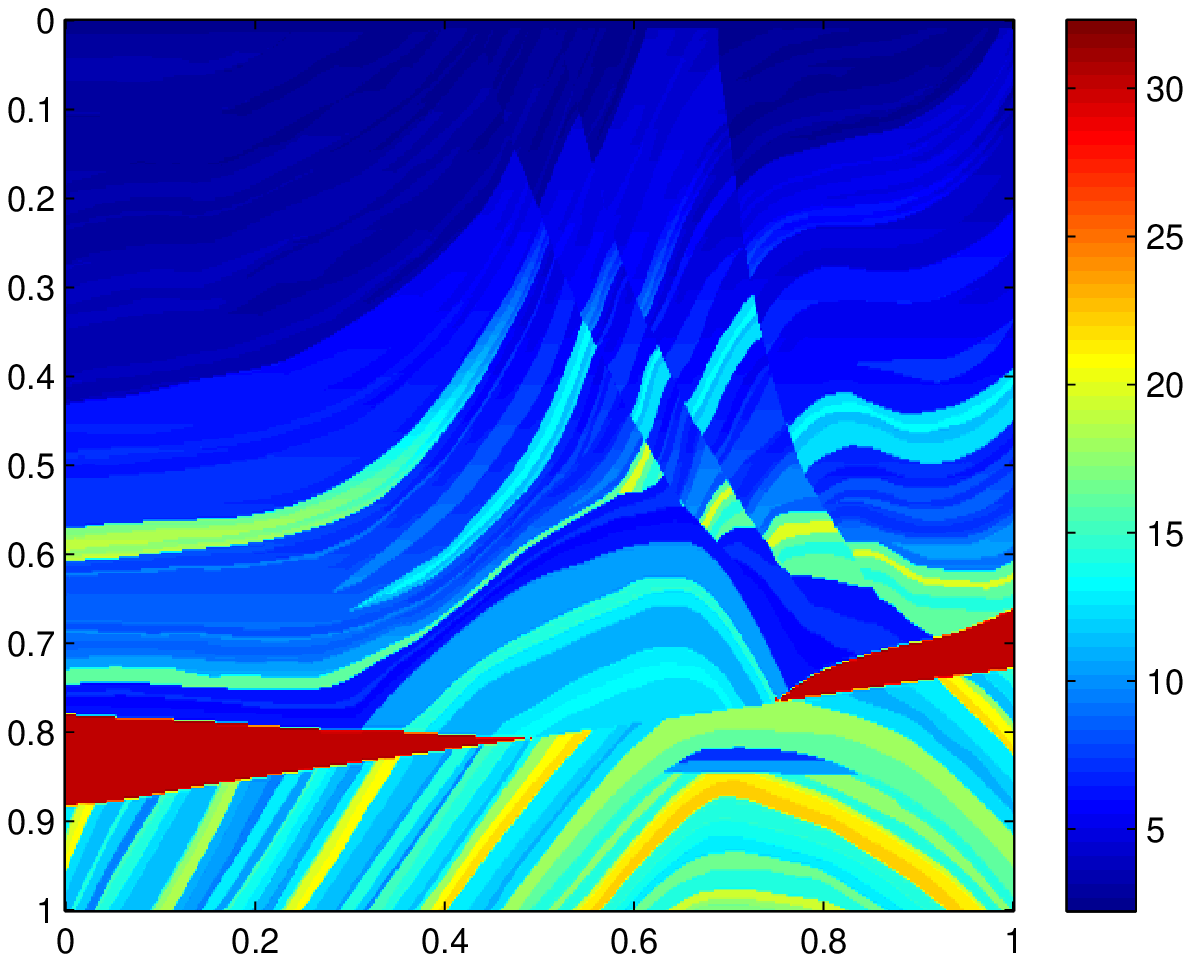} \includegraphics[scale=0.5]{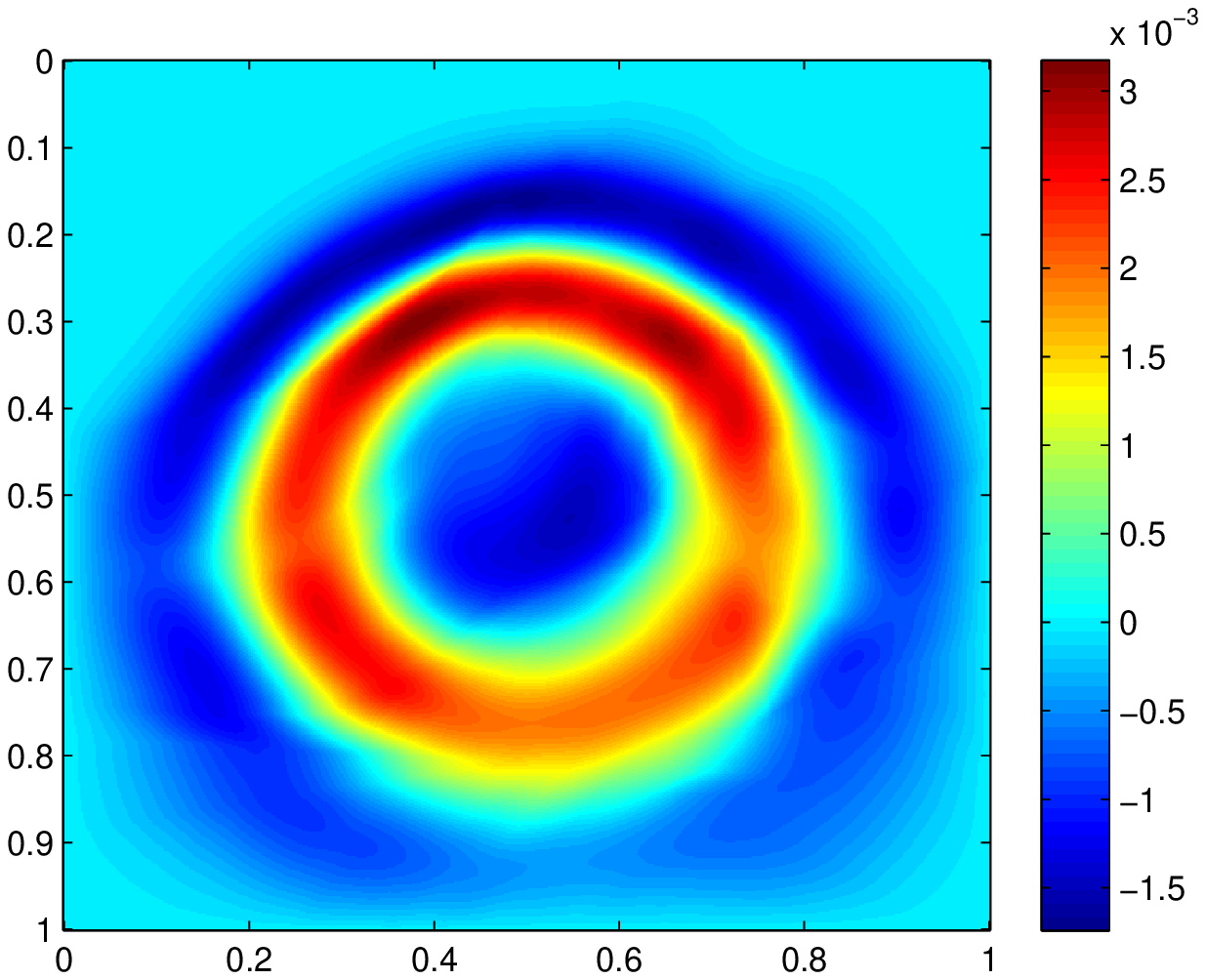}

\caption{Left: a subset of the Marmousi model. Right: fine grid solution.}

\label{fig:marmousi}
\end{figure}

In this section, we will present some numerical examples to show the
performance of our multiscale method. The media that we will consider is a
 heterogeneous field which is a modified
 Marmousi model (see the left plot of Figure \ref{fig:marmousi}).
We have also considered more regular periodic
highly heterogeneous fields and observed similar results.
 We will compare both
the accuracy and efficiency of our method with the
direct fine scale simulation defined in (\ref{eq:globalFEM}).
To compare the accuracy, we will use the following error quantities
\[
e_{2}=\cfrac{\|u_{H}-u\|_{L^{2}(\Omega)}}{\|u\|_{L^{2}(\Omega)}},\quad\overline{e_{2}}=\cfrac{\sqrt{\sum_{K\in\mathcal{T}^{H}}|\int_{K}u_{H}-\int_{K}u|^{2}}}{\sqrt{\sum_{K\in\mathcal{T}^{H}}|\int_{K}u|^{2}}},\quad e_{H^{1}}=\cfrac{\|\nabla(u_{H}-u)\|_{L^{2}(\Omega)}}{\|\nabla u\|_{L^{2}(\Omega)}}\]
 which are the relative $L^{2}$ norm error, the relative $L^{2}$
norm error for coarse grid averages and the relative $L^{2}$ norm
error of the gradient. We will also consider the jump error on coarse
grid edges defined by \[
e_{Jump}=\sum_{e\in\mathcal{E}^{H}}\int_{e}[u]_{e}^{2}.\]

Moreover, we let $t_{off}$ be the time needed for offline computations
and $t_{on}$ be the online computational time. These quantities are
used to compare the efficiency of our method with direct fine scale
simulation. To perform a fair comparison, we will use the same time
step size for both of our GMsFEM and the fine scale method, since
we only consider spatial upscaling in this paper.
However, we note that multiscale basis functions can be used
for different source terms and boundary conditions which will
provide a substantial computational saving.
 Furthermore, we
will take $\gamma=2$ and $\Omega=[0,1]^{2}$ for all of our examples.
The initial conditions $g_{0}$ and $g_{1}$ are zero. Throughout
the paper, all computational times are measured in seconds.

%We consider wave propagation in the medium defined
%by a subset of the Marmousi model, where the coefficient $a(x,y)$
%is shown in the left plot of Figure \ref{fig:marmousi}.
The Ricker wavelet with frequency $f_{0}=20$
\begin{align*}
f(x,y) & =(10)^{2}e^{-10^{2}((x-0.5)^{2}+(y-0.5)^{2})}(1-2\pi^{2}f_{0}^{2}(t-2/f_{0})^{2})e^{-\pi^{2}f_{0}^{2}(t-2/f_{o})^{2}}\end{align*}
is used as the source term.
We will compute the solution at time $T=0.4$.
The coarse mesh size is taken as $H=1/16$. Each coarse
grid block is divided into a $32\times32$ grid, that is, $n=32$.
Thus, the fine mesh size $h=1/512$ and there are totally $128$ and
$961$ local basis functions in the space $V_{H}^{1}(K)$
and $V_{H}^{2}(K)$ respectively on each coarse grid block.
The time step size for both GMsFEM and the fine grid solver
is taken as $\Delta t=h/80$ in order to meet the stability requirement
and the computation time for fine grid solution is $55.06$. We will
compare the accuracy and efficiency of our method using the solution
computed at the time $T=0.2$, which is shown in
the right figure of Figure
\ref{fig:marmousi}.

%\begin{figure}[H]
%\centering \includegraphics[scale=0.5]{sol}

%\caption{Exact solution.}

%\label{fig:sol} %
%\end{figure}

%For the numerical computations by our multiscale method, we take $\Delta t=h/4$
%in order to meet the stability requirement and the other parameters
%remain the same.

\begin{figure}[H]
 \centering

\includegraphics[scale=0.6]{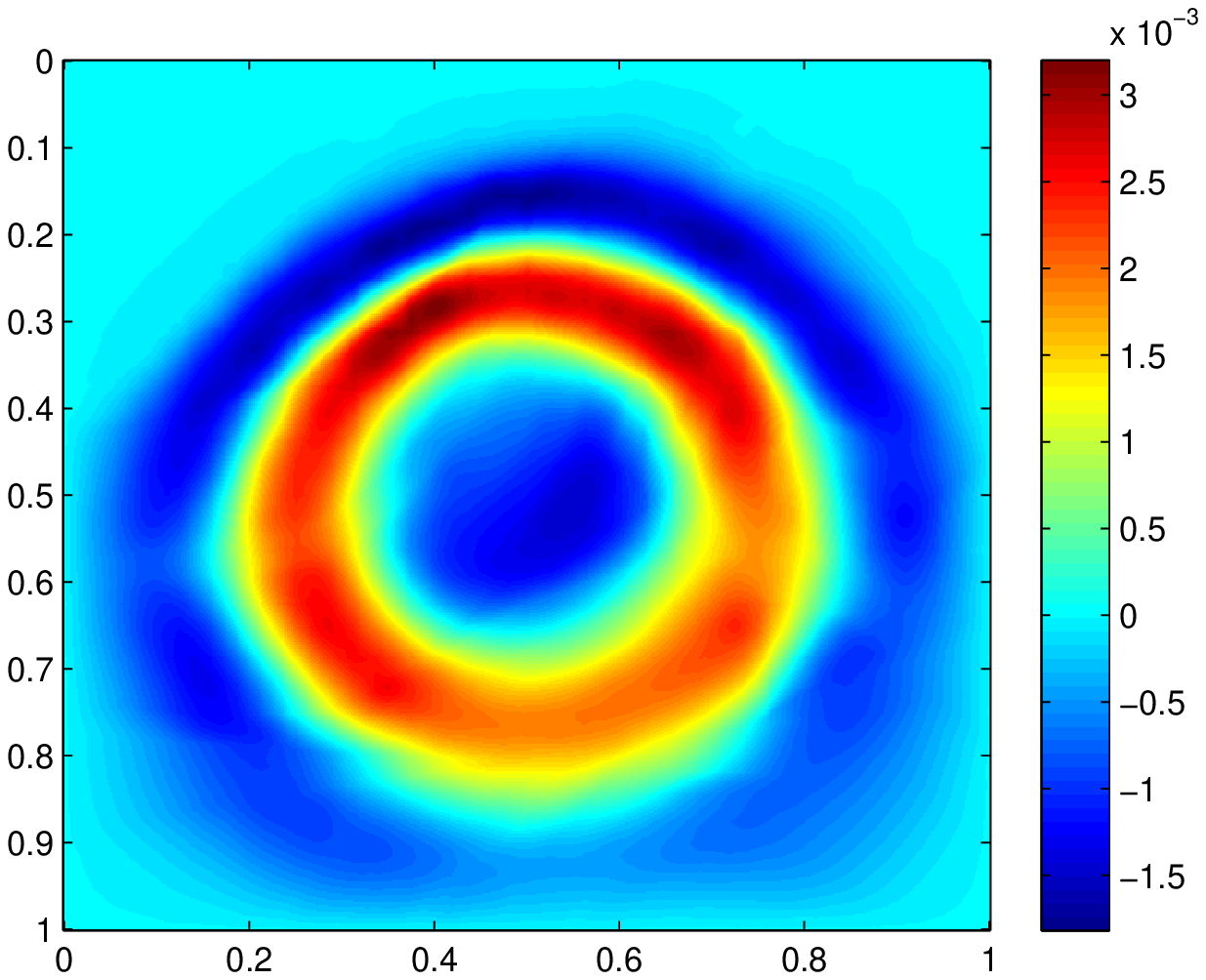}\includegraphics[scale=0.6]{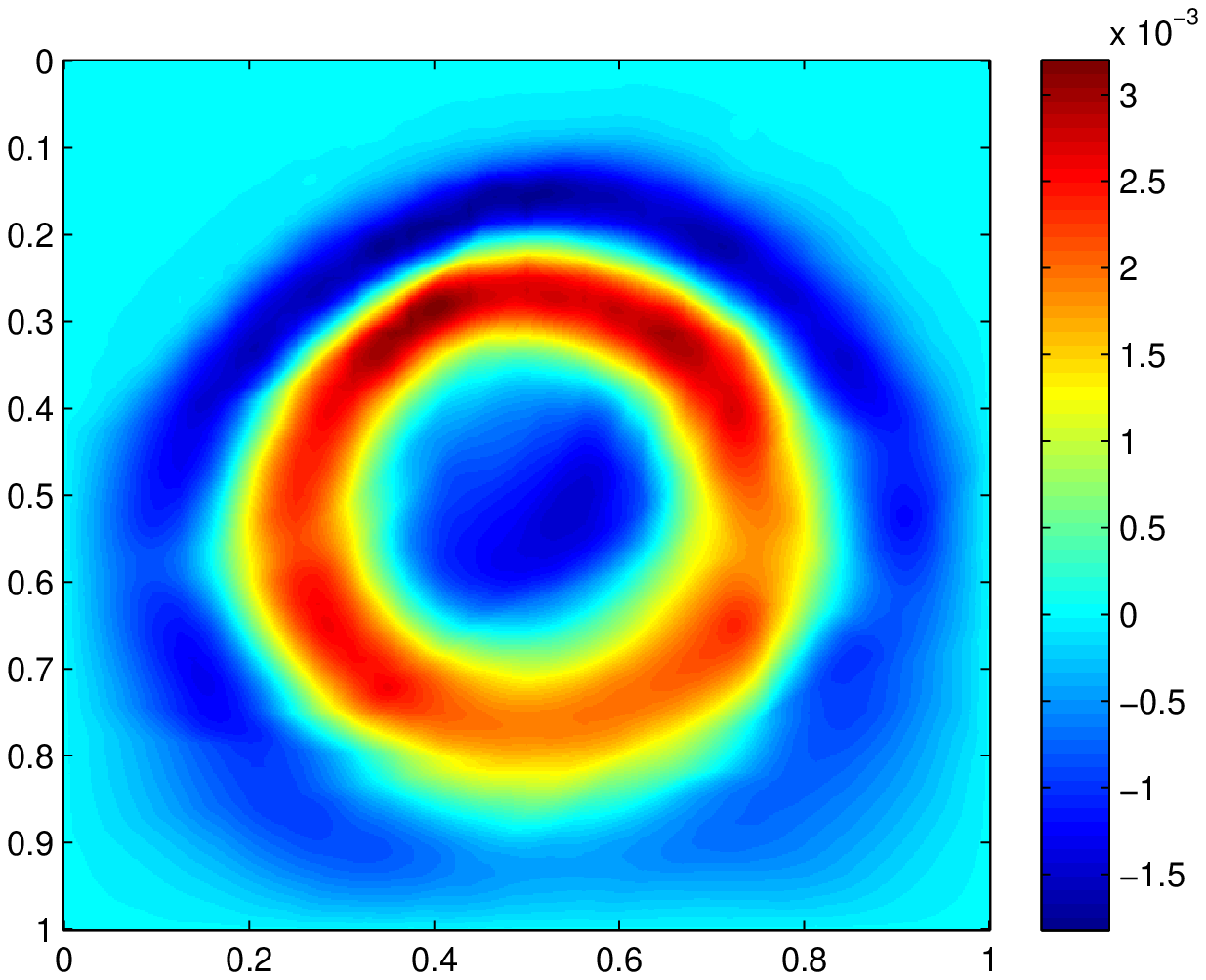}
\caption{Left: $75\%$ energy. Right: $80\%$ energy.}

\label{fig:ex3}
\end{figure}

\begin{figure}[H]
 \centering{}\includegraphics[scale=0.4]{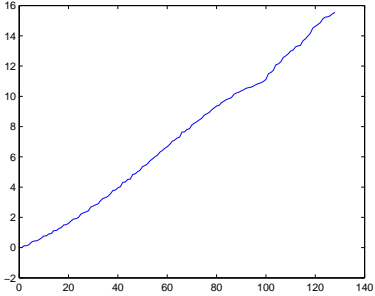} \caption{Eigenvalues for the space $V_{H}^{1}$.}

\label{fig:eig3}
\end{figure}

In Table \ref{tab:ex3} and Table \ref{tab:ex3new}, we present the
errors and computational times for the case with $m=1$, that is,
we only use the first eigenfunction in the space $V_{H}^{2}$.
%We
%see that, if we use $85\%$ of the total energy for the space of harmonic
%extensions $V_{H}^{1}$, then we need $45$ to $54$ basis functions
%on each coarse grid. The computational time for the offline procedures
%is $3533.27$ and the time for online computations is $55.24$. For
%this case, we see that the computational time is about the same as
%that of the direct fine grid simulation. When
We see that if we use
$80\%$ of the total
energy, the number of basis functions is between
$33$ and $40$ on each coarse grid
while the computational time for the offline procedure
is $1019.06$ and the time for online computations is
$32.43$. Note that the online computational time is about $59\%$
of that of the online computational time of the direct fine grid simulation.
%Moreover, in this case, t
The relative $L^{2}$ error and the relative
error for cell averages are only $3.92\%$ and $2.74\%$ respectively.
In addition, the relative error for the gradient is $14.86\%$ and
the jump error is $0.003$. When $75\%$ of the total energy is
used, the number of basis functions is  reduced to a number
between $24$ and $29$ while the computational time for the offline
procedures is $326.83$ and the time for online computations is $18.21$.
The time for the online computation is $33\%$ of the
time required for direct fine grid simulation. The relative $L^{2}$
error and the relative error for cell averages are increased slightly
to $4.23\%$ and $3.12\%$, respectively.
In Table \ref{tab:ex3}, we also present the values of $\mu_{min}$
for the space $V_{H}^{1}$. Moreover, the eigenvalues are shown in
Figure \ref{fig:eig3}.
The numerical solutions for
these  cases are shown in Figure \ref{fig:ex3}.
%Comparing the
%fine scale solution shown in Figure \ref{fig:marmousi}, we see that
%our method is able to capture the solution well.
We note that the error decay is not fast mostly due to
the error contribution because of the modes corresponding to
the interior. Even though the error between the GMsFEM solution
and the solution computed using the entire snapshot space
$V_H^1$ is very small, the overall error between the GMsFEM solution
and the fine-scale solution may not be small because we have only used
one basis function in $V_H^2$. Next, we will add more basis functions
from $V_H^2$ and compare the errors.

\begin{table}[H]
 \centering \begin{tabular}{|c|c|c|c|c|c|c|}
\hline
Energy  & Number of basis  & $e_{2}$  & $\overline{e_{2}}$  & $e_{H^{1}}$  & $e_{Jump}$  & $\mu_{min}$\tabularnewline
\hline
$75\%$  & 24-29  & 0.0423  & 0.0312  & 0.1542  & 4.7304e-04  & 1.9414\tabularnewline
\hline
$80\%$  & 33-40  & 0.0392  & 0.0274  & 0.1486  & 3.0671e-04  & 2.9992\tabularnewline
\hline
\hline
\end{tabular}\caption{Errors for various choices of energy for the space $V_{H}^{1}$.}

\label{tab:ex3}
\end{table}

\begin{table}[H]
 \centering \begin{tabular}{|c|c|c|}
\hline
Energy  & $t_{off}$  & $t_{on}$ \tabularnewline
\hline
$75\%$  & 326.83  & 18.21\tabularnewline
\hline
$80\%$  & 1019.06  & 32.43\tabularnewline
\hline
\end{tabular}\caption{Offline and online computational times.}

\label{tab:ex3new}
\end{table}

%When $85\%$ of the total energy is used, the eigenvalue is $4.44$,
%when $80\%$ of the total energy is used, the eigenvalue is $3.00$
%and when $75\%$ of the total energy is used, the eigenvalue is $1.94$.

%

\begin{figure}[H]
 \centering \includegraphics[scale=0.4]{numesol_case3_m1_E0025}\includegraphics[scale=0.4]{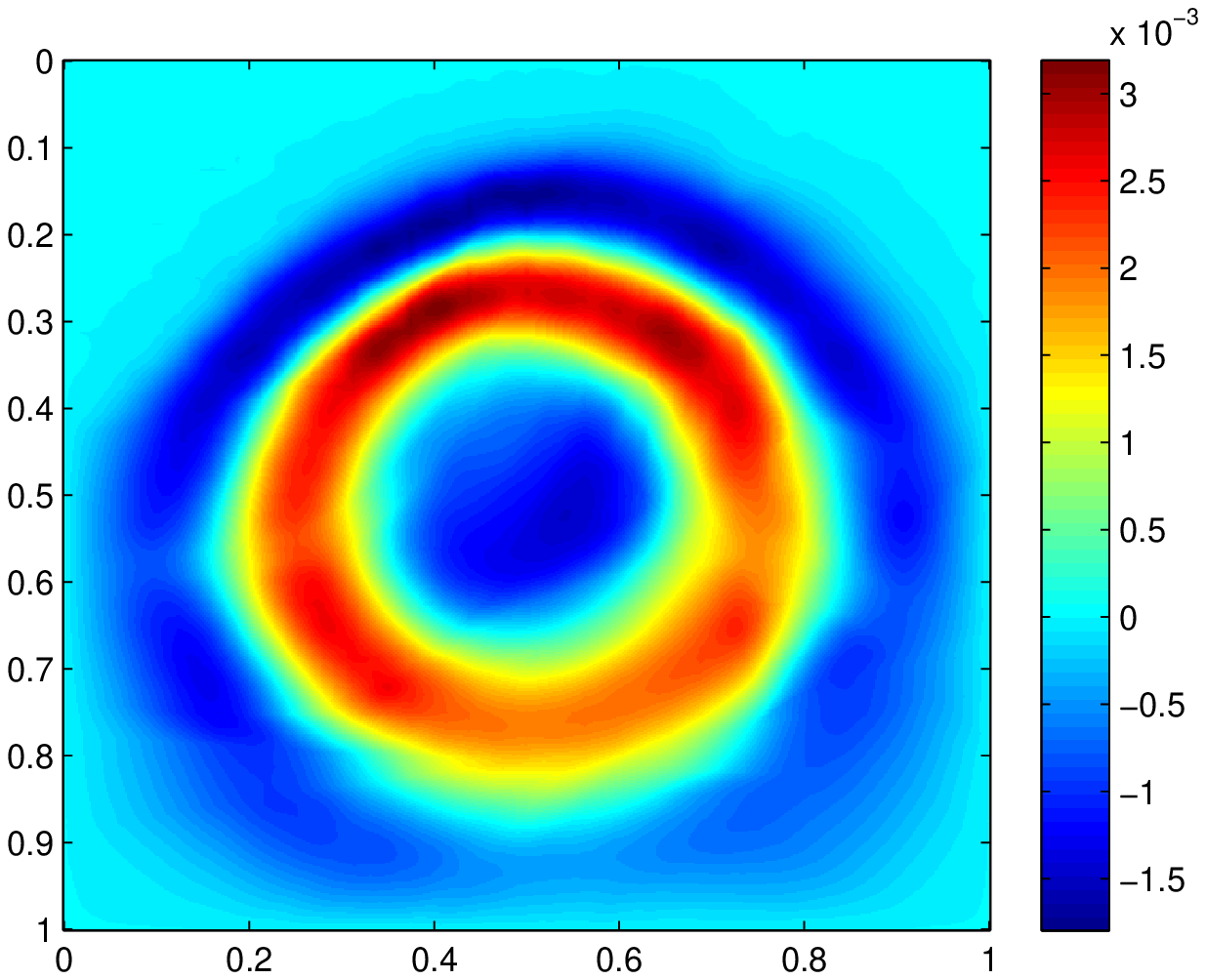}
\includegraphics[scale=0.4]{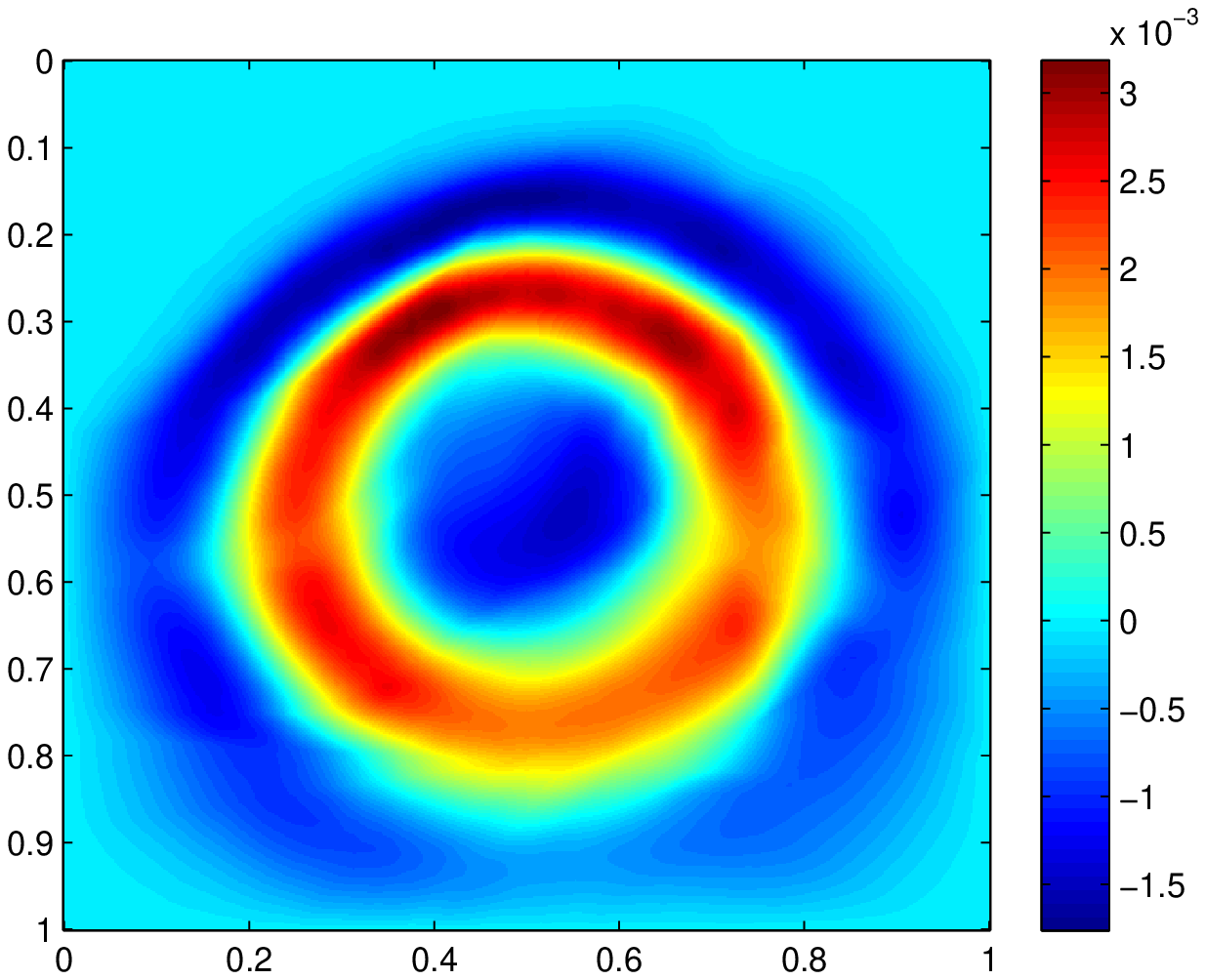}\includegraphics[scale=0.4]{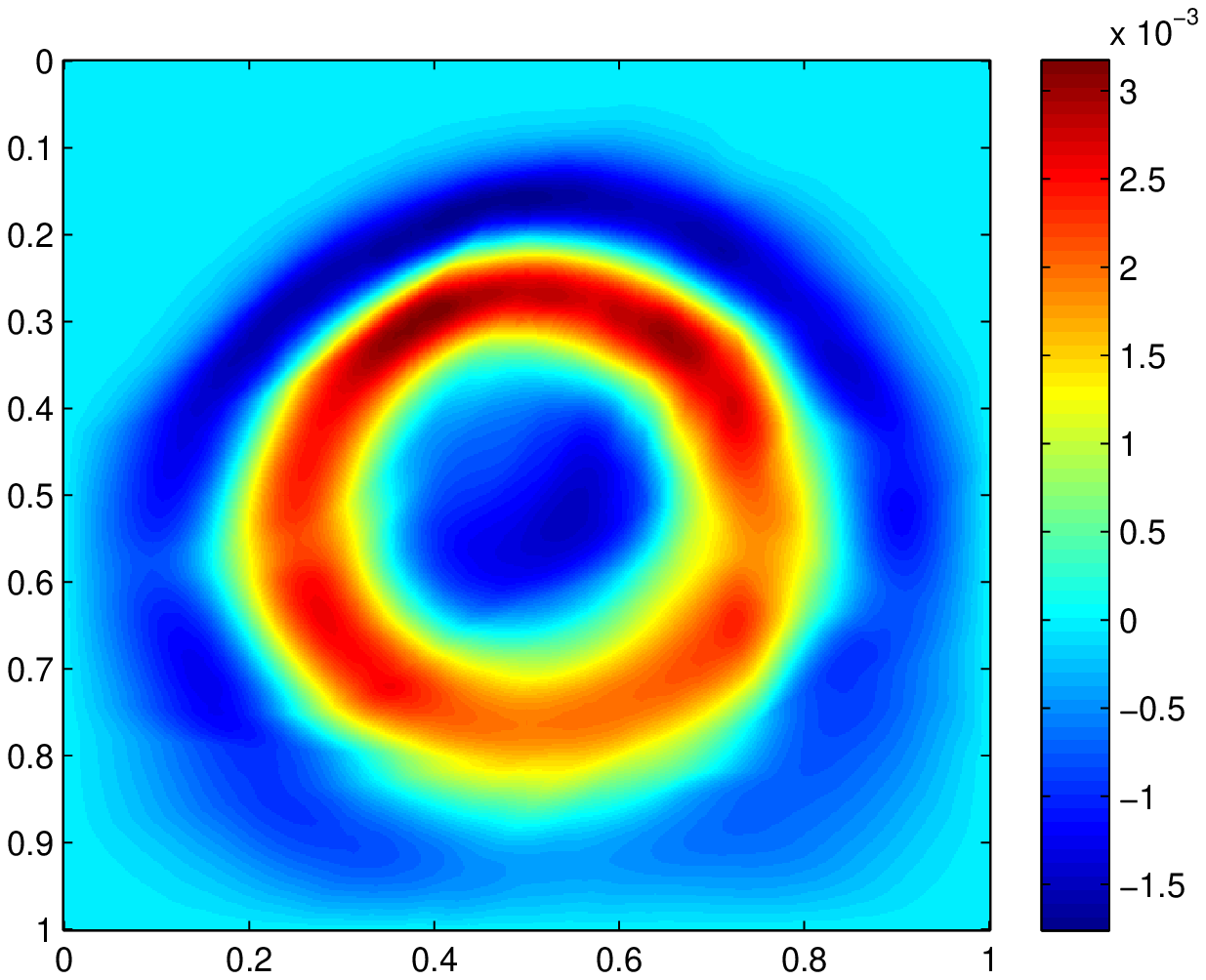}

\caption{Numerical solutions with various number of eigenfunctions $(m)$ in
the space $V_{H}^{2}$ for using $75\%$ energy in $V_{H}^{1}$.
Upper-left: $m=1$. Upper-right: $m=2$. Lower-left: $m=3$. Lower-right:
$m=5$.}

\label{fig:ex3a}
\end{figure}

Next, we will investigate the use of more eigenfunctions in the space
$V_{H}^{2}$ that will allow reducing the overall error.
 To do so, we consider the first case where $75\%$ energy in the space
$V_{H}^{1}$ is used
 and we consider using various number of eigenfunctions in
$V_{H}^{2}$. The errors and computational times are shown in Table
\ref{tab:ex3a} and Table \ref{tab:ex3anew}. In general, we obtain
better numerical approximations as more eigenfunctions are used. When
two eigenfunctions are used (this corresponds to using less than
$3\%$ of the total local degrees of freedom
 in constructing all GMsFEM basis functions),
the relative error is $3.52\%$ and the
online computational time is $18.64$. When five eigenfunctions are
used, the relative error is $1.93\%$ and the online computational
time is $18.21$. Thus, we see that adding a few eigenfunctions in
the space $V_{H}^{2}$ will improve the multiscale solution.
This indicates that for the multiscale wave simulations, the modes
that represent the interior nodes can improve the accuracy of
the method and play an important role in obtaining an accurate solution.
The numerical solution for these $4$ cases are shown in Figure \ref{fig:ex3a}.
We see that our method is able to capture the solution well. We also
report the largest eigenvalue used in Table \ref{tab:ex3a}.

\begin{table}[H]
 \centering \begin{tabular}{|c|c|c|c|c|c|}
\hline
m  & $e_{2}$  & $\overline{e_{2}}$  & $e_{H^{1}}$  & $e_{Jump}$  & $\lambda_{min}$\tabularnewline
\hline
1  & 0.0423  & 0.0312  & 0.1542  & 4.7304e-04  & 3.4805e+04\tabularnewline
\hline
2  & 0.0352  & 0.0259  & 0.1346  & 4.7030e-04  & 3.4873e+04\tabularnewline
\hline
3  & 0.0227  & 0.0187  & 0.0945  & 4.5931e-04  & 5.5906e+04\tabularnewline
\hline
5  & 0.0193  & 0.0163  & 0.0833  & 4.5910e-04  & 6.9650e+04\tabularnewline
\hline
\end{tabular}\caption{Errors for various number of eigenfunctions in $V_{H}^{2}$ for using
$75\%$ energy in $V_{H}^{1}$.}

\label{tab:ex3a}
\end{table}

\begin{table}[H]
 \centering \begin{tabular}{|c|c|c|}
\hline
m  & $t_{off}$  & $t_{on}$ \tabularnewline
\hline
1  & 326.83  & 18.21\tabularnewline
\hline
2  & 368.89  & 18.64\tabularnewline
\hline
3  & 405.73  & 19.88\tabularnewline
\hline
5  & 528.47  & 25.96\tabularnewline
\hline
\end{tabular}\caption{Offline and online computational times for various number of eigenfunctions
in $V_{H}^{2}$ for using $75\%$ energy in $V_{H}^{1}$.}

\label{tab:ex3anew}
\end{table}

We would like to remark that the computational gain will be higher when
implicit methods are used or we employ finer grids to resolve the problem.
In the latter case, the CPU time for coarse-grid simulations will not
change.

\subsection{The use of oversampling}

In this section, we present the performance of the method
when the basis functions in the space $V_H^1$ are obtained by oversampling.
We consider the previous example.
The oversampling technique is used and
the harmonic extension problems are solved on enlarged coarse grids, which are obtained
by extending the original coarse grids by $H/16$ on each side.
The results for using one basis functions in $V_H^2$
and various number of basis functions in $V_H^1$ are shown in Table \ref{tab:over_ex2}.
Moreover, we compute the errors using $73\%$ energy for $V_H^1$ and various number of basis functions in $V_H^2$.
The results are presented in Table \ref{tab:over_ex2a}.
We observe that there is no improvement in this case. This is due
to the error  from the
 modes representing internal nodes.

\begin{table}[ht]
\centering \begin{tabular}{|c|c|c|c|c|c|c|}
\hline
Energy  & Number of basis  & $e_{2}$  & $\overline{e_{2}}$  & $e_{H^{1}}$  & $e_{Jump}$  & $\mu_{min}$\tabularnewline
\hline
$73\%$  & 24-30 & 0.0673 & 0.0583 & 0.1866 & 5.2038e-04 & 1.6755\tabularnewline
\hline
$79\%$  & 33-40 & 0.0640 & 0.0548 & 0.1827 & 3.4797e-04 & 2.5681\tabularnewline
\hline
$84\%$  & 45-55 & 0.0626 & 0.0534 & 0.1809 & 2.6388e-04 & 3.7918\tabularnewline
\hline
\end{tabular}\caption{Simulation results with one basis function in $V_H^2$.}

\label{tab:over_ex2} %
\end{table}

\begin{table}[ht]
\centering \begin{tabular}{|c|c|c|c|c|c|}
\hline
m  & $e_{2}$  & $\overline{e_{2}}$  & $e_{H^{1}}$  & $e_{Jump}$  & $\lambda_{min}$\tabularnewline
\hline
1  & 0.0673  & 0.0583  & 0.1866  & 5.2038e-04  & 3.4805e+04\tabularnewline
\hline
2  & 0.0596  & 0.0524  & 0.1666  & 5.1865e-04  & 3.4873e+04\tabularnewline
\hline
3  & 0.0488  & 0.0449  & 0.1332  & 5.0929e-04  & 5.5906e+04\tabularnewline
\hline
5  & 0.0449  & 0.0419  & 0.1220  & 5.0793e-04  & 6.9650e+04\tabularnewline
\hline
\end{tabular}\caption{Errors and computational times for various number of eigenfunctions
in $V_{H}^{2}$ for using $E=73\%$.}

\label{tab:over_ex2a}%
\end{table}

\section{Stability and convergence}

\label{sec:stab}

In this section, we will prove the stability and convergence of the
generalized multiscale finite element method constructed in Section
\ref{sec:method}. We will first state and prove some preliminary
results, and then prove the main convergence theorem for the semi-discrete
scheme (\ref{eq:ipdg}).

\subsection{Preliminaries}

Before we analyze the convergence of our GMsFEM,
we first prove some basic results. To do so, we introduce some notations
and state the assumptions required in our analysis. For functions
$u,v\in H^{1}(\mathcal{T}^{H})$, we define the bilinear form $a(\cdot,\cdot)$
by \[
a(u,v)=\sum_{K\in\mathcal{T}^{H}}\int_{K}a\nabla u\cdot\nabla v.\]
 Moreover, for any function $u\in H^{1}(\mathcal{T}^{H})$, we define
the $a$-norm by \[
\|u\|_{a}=\left(a(u,u)+\cfrac{\gamma}{h}\sum_{e\in\mathcal{E}^{H}}\|a^{\frac{1}{2}}[u]_{e}\|_{L^{2}(e)}^{2}\right)^{\frac{1}{2}}\]
 and the $a$-semi-norm by \[
|u|_{a}=a(u,u)^{\frac{1}{2}}.\]
 Furthermore, the broken $H^{1}$-norm for $u\in H^{1}(\mathcal{T}^{H})$
is defined as \[
\|u\|_{H^{1}(\mathcal{T}^{H})}=\left(\sum_{K\in T^{H}}|u|_{H^{1}(K)}^{2}+\cfrac{\gamma}{h}\sum_{e\in\mathcal{E}^{H}}\|[u]_{e}\|_{L^{2}(e)}^{2}\right)^{\frac{1}{2}}.\]

\begin{assumption} The function $a(x)$ is bounded, that is, there
exist positive numbers $a_{0}$ and $a_{1}$ such that \[
a_{0}\leq a(x)\leq a_{1},\quad\quad\forall\, x\in\Omega.\]
 This assumption implies that the norms $\|\cdot\|_{a}$ and $\|\cdot\|_{H^{1}(\mathcal{T}^{H})}$
are equivalent. \end{assumption}

In the following, we will describe the consistency of the method (\ref{eq:ipdg}).
We define the consistency error by \begin{equation}
R_{u_{h}}(v)=l(v)-(\frac{\partial^{2}u_{h}}{\partial t^{2}},v)-a_{DG}(u_{h},v),\quad\quad\forall\, v\in V_{H},\end{equation}
 where $u_{h}$ is the fine grid finite element solution defined in (\ref{eq:globalFEM}).
Clearly, we have \begin{equation}
R_{u_{h}}(v)=0,\quad\quad\forall\, v\in V_{H}^{2}\end{equation}
since $V_H^2 \subset V_h$.
 Thus, we only need to estimate $R_{u_{h}}(v)$ for $v\in V_{H}^{1}$.
The following lemma states that the method (\ref{eq:ipdg}) is consistent
with the fine grid solution defined by (\ref{eq:globalFEM}). The
proof will be presented in the Appendix.

\begin{lemma} \label{lem:consistent} Let $u_{h}$ and $u$ be the
finite element solution defined in (\ref{eq:globalFEM}) and the exact
solution of the wave propagation problem (\ref{eq:waveeqn}) respectively.
%\begin{equation}
%\sum_{K\in\mathcal{T}^{H}}(\frac{\partial u}{\partial t^{2}},v)_{L^{2}(K)}+a_{DG}(u,v)=\sum_{K\in\mathcal{T}^{H}}(F,v)_{L^{2}(K)},\quad\forall\, v\in V^1_{H}\label{eq:consistency}%\end{equation}
If $u\in H^{2}(\Omega)$, then we have \begin{equation}
|R_{u_{h}}(v)|\leq C(u,f)h\|v\|_{a},\quad\quad v\in V_{H}^{1}\label{eq:consistency}\end{equation}
 where $C(u,f)$ is a constant which depends on the solution $u$
and the source term $f$ but independent of the fine mesh size $h$.
This inequality gives the consistency of our method. \end{lemma}

Next we will prove that the bilinear form $a_{DG}$ satisfies the
following coercivity and continuity conditions for suitably chosen
penalty parameter $\gamma>0$.

\begin{lemma} Let $\gamma$ be sufficiently large. Then we have \begin{equation}
\frac{1}{2}\,\|v\|_{a}^{2}\leq a_{DG}(v,v),\quad\forall\, v\in V_{H}\label{eq:coercivity}\end{equation}
 and \begin{equation}
a_{DG}(u,v)\leq2\|u\|_{a}\|v\|_{a},\quad\forall\, u,v\in V_{H}.\label{eq:continuity (bilinear)}\end{equation}
 \end{lemma}

\noindent \textit{Proof: } By the definition of $a_{DG}$, we have
\[
a_{DG}(v,v)=\sum_{K\in\mathcal{T}^{H}}\int_{K}a\nabla v\cdot\nabla v+\sum_{e\in\mathcal{E}^{H}}\Big(-2\int_{e}\{a\nabla v\cdot n\}_{e}\,[v]_{e}+\cfrac{\gamma}{h}\int_{e}a[v]_{e}^{2}\Big), \]
and by the definition of the $a$-norm, we have
\begin{equation}
a_{DG}(v,v)=\|v\|_{a}^{2}-2\sum_{e\in\mathcal{E}^{H}}\Big(\int_{e}\{a\nabla v\cdot n\}_{e}\,[v]_{e}\Big). \label{eq:coer_1}
\end{equation}
 Using the Cauchy-Schwarz inequality, we obtain \[
2\sum_{e\in\mathcal{E}^{H}}\Big(\int_{e}\{a\nabla v\cdot n\}_{e}\,[v]_{e}\Big)\leq\frac{2h}{\gamma a_{0}}\sum_{e\in\mathcal{E}^{H}}\int_{e}\{a\nabla v\cdot n\}_{e}^{2}+\frac{a_{0}}{2}\sum_{e\in\mathcal{E}^{H}}\frac{\gamma}{h}\int_{e}[v]_{e}^{2}\]
 which implies \[
2\sum_{e\in\mathcal{E}^{H}}\Big(\int_{e}\{a\nabla v\cdot n\}_{e}\,[v]_{e}\Big)\leq\frac{h}{\gamma a_{0}}\sum_{K\in\mathcal{T}^{H}}\int_{\partial K}(a\nabla v\cdot n)^{2}+\frac{1}{2}\sum_{e\in\mathcal{E}^{H}}\frac{\gamma}{h}\int_{e}a[v]_{e}^{2}.\]
 Since $v$ is a piecewise linear function, there is a uniform constant $\Lambda>0$ such that
 \begin{equation}
2h\sum_{K\in\mathcal{T}^{H}}\int_{\partial K}(a\nabla v\cdot n)^{2}\leq \Lambda a_{1}|v|_{a}^{2}.\label{eq:discretetrace}\end{equation}
 So we have \[
2\sum_{e\in\mathcal{E}^{H}}\Big(\int_{e}\{a\nabla v\cdot n\}_{e}\,[v]_{e}\Big)\leq\frac{\Lambda a_{1}}{2\gamma a_{0}}|v|_{a}^{2}+\frac{1}{2}\sum_{e\in\mathcal{E}^{H}}\frac{\gamma}{h}\int_{e}a[v]_{e}^{2}.\]
 Therefore, from (\ref{eq:coer_1}), we obtain \[
\frac{1}{2}\,\|v\|_{a}^{2}\leq a_{DG}(v,v),\quad\forall\, v\in V_{H}\]
 if we take $\gamma\geq \Lambda a_{1}a_{0}^{-1}$. Thus, we have proved (\ref{eq:coercivity}).

We can prove (\ref{eq:continuity (bilinear)}) in the similar way.
By the definition of $a_{DG}$, we have \begin{eqnarray*}
 &  & \left|a_{DG}(u,v)\right|\\
 & = & \left|\sum_{K\in\mathcal{T}^{H}}\int_{K}a\nabla u\cdot\nabla v+\sum_{e\in\mathcal{E}^{H}}\Big(-\int_{e}\{a\nabla u\cdot n\}_{e}\cdot[v]_{e}-\int_{e}\{a\nabla v\cdot n\}_{e}\cdot[u]_{e}+\cfrac{\gamma}{h}\int_{e}a[u]_{e}\cdot[v]_{e}\Big)\right|\\
 & \leq & I_{1}+I_{2}+I_{3}+I_{4},\end{eqnarray*}
 where \begin{align*}
I_{1} & =\sum_{K\in\mathcal{T}^{H}}|\int_{K}a\nabla u\cdot\nabla v|,\quad\quad I_{2}=\sum_{e\in\mathcal{E}^{H}}|\int_{e}\{a\nabla u\cdot n\}_{e}\cdot[v]_{e}|,\\
I_{3} & =\sum_{e\in\mathcal{E}^{H}}|\int_{e}\{a\nabla v\cdot n\}_{e}\cdot[u]_{e}|,\quad\quad I_{4}=\sum_{e\in\mathcal{E}^{H}}\cfrac{\gamma}{h}\int_{e}|a[u]_{e}\cdot[v]_{e}|.\end{align*}
 First, we note that $I_{1}$ and $I_{4}$ can be estimated easily
as follows: \begin{align*}
I_{1} & \leq\left(\sum_{K\in\mathcal{T}^{H}}\int_{K}a|\nabla u|^{2}\right)^{\frac{1}{2}}\left(\sum_{K\in\mathcal{T}^{H}}\int_{K}a|\nabla v|^{2}\right)^{\frac{1}{2}},\\
I_{4} & \leq\left(\sum_{e\in\mathcal{E}^{H}}\cfrac{\gamma}{h}\int_{e}a[u]_{e}^{2}\right)^{\frac{1}{2}}\left(\sum_{e\in\mathcal{E}^{H}}\cfrac{\gamma}{h}\int_{e}a[v]_{e}^{2}\right)^{\frac{1}{2}}.\end{align*}
 For $I_{2}$, we can estimate as follows \begin{align*}
I_{2} & \leq\left(\cfrac{h}{\gamma}\sum_{e\in\mathcal{E}^{H}}\int_{e}\{a\nabla u\cdot n\}_{e}^{2}\right)^{\frac{1}{2}}\left(\sum_{e\in\mathcal{E}^{H}}\cfrac{\gamma}{h}\int_{e}[v]_{e}^{2}\right)^{\frac{1}{2}}\\
 & \leq\left(\cfrac{h}{\gamma a_{0}}\sum_{K\in\mathcal{T}^{H}}\int_{\partial K}(a\nabla u\cdot n)^{2}\right)^{\frac{1}{2}}\left(\sum_{e\in\mathcal{E}^{H}}\cfrac{\gamma}{h}\int_{e}a[v]_{e}^{2}\right)^{\frac{1}{2}}\\
 & \leq(\frac{\Lambda a_{1}}{\gamma a_{0}})^{\frac{1}{2}}|u|_{a}\left(\sum_{e\in\mathcal{E}^{H}}\cfrac{\gamma}{h}\int_{e}[v]_{e}^{2}\right)^{\frac{1}{2}}.\end{align*}
 The same idea can be applied to estimate $I_{3}$ to obtain \begin{align*}
I_{3} & \leq(\frac{\Lambda a_{1}}{\gamma a_{0}})^{\frac{1}{2}}|v|_{a}\left(\sum_{e\in\mathcal{E}^{H}}\cfrac{\gamma}{h}\int_{e}a[u]_{e}^{2}\right)^{\frac{1}{2}}.\end{align*}
 Finally, combining the above estimates, we have \[
\left|a_{DG}(u,v)\right|\leq2\|u\|_{a}\|v\|_{a}\]
 when $\gamma\geq \Lambda a_{1}a_{0}^{-1}.$

\begin{flushright}
$\square$
\par\end{flushright}

Next, we will prove the convergence of the semi-discrete scheme (\ref{eq:ipdg}).
First, we define the following error quantities. Let
\begin{equation}
\eta=u_{h}-w_{H},\quad\xi=u_{H}-w_{H},\quad\text{and}\quad \varepsilon=u_{h}-u_{H},
\label{eq:errordef}
\end{equation}
 where $w_{H}\in V_{H}$ is defined by solving the following elliptic
projection problem \begin{equation}
a_{DG}(w_{H},v)=a_{DG}(u_{h},v)+R_{u_{h}}(v),\quad\quad\forall\, v\in V_{H}.\label{eq:ellipticproj}\end{equation}
Notice that $\varepsilon$ is the difference between the multiscale solution $u_H$
and the fine grid finite element solution $u_h$.
Moreover, $\eta$ measures the difference between the fine grid solution $u_h$ as its projection $w_H$.
In the following, we will prove estimates for $\varepsilon$.
First, we let
\[
\|\varepsilon\|_{L^{\infty}([0,T];L^{2}(\Omega))}=\max_{0\leq t\leq T}\|\varepsilon\|_{L^{2}(\Omega)}
\quad\text{and}\quad\|\varepsilon\|_{L^{\infty}([0,T];a)}=\max_{0\leq t\leq T}\|\varepsilon\|_{a}.\]
Then we will prove the following two inequalities,
which estimate the error for the solution $\varepsilon$ by the error for the projection $\eta$
and the initial errors $\mathcal{I}_1$ and $\mathcal{I}_2$,
which are defined in the statements of the theorems.

\begin{theorem} Let $\varepsilon,\eta$ and $\xi$ be the error quantities
defined in (\ref{eq:errordef}). Then we have the following error
bound \begin{equation}
\begin{split} & \:\|\varepsilon_{t}\|_{L^{\infty}([0,T];L^{2}(\Omega))}+\|\varepsilon\|_{L^{\infty}([0,T];a)}\\
\leq & \: C\Big(\|\eta_{t}\|_{L^{\infty}([0,T];L^{2}(\Omega))}+\|\eta\|_{L^{\infty}([0,T];H^{1}(\mathcal{T}^{H}))}+\|\eta_{tt}\|_{L^{1}([0,T];L^{2}(\Omega))} + \mathcal{I}_1 \Big),
\end{split}
\label{eq:errorineq1}\end{equation}
where $\mathcal{I}_1 = \|\xi_{t}(0)\|_{L^{2}(\Omega)}+\|\xi(0)\|_{H^{1}(\mathcal{T}^{H})}$.
 \label{thm:est1} \end{theorem} %
%\marginpar{Should we cross reference to $\xi$ and $\eta$}
\textit{Proof:} First, using (\ref{eq:ipdg}) and the definition
of $\xi$, we have \[
(\xi_{tt},v)+a_{DG}(\xi,v)=(f,v)-((w_{H})_{tt},v)-a_{DG}(w_{H},v).\]
 Then by (\ref{eq:consistency}), we have \begin{equation}
(\xi_{tt},v)+a_{DG}(\xi,v)=(\eta_{tt},v).\label{eq:1}\end{equation}
 Taking $v=\xi_{t}$ in (\ref{eq:1}), we have \begin{align*}
(\xi_{tt},\xi_{t})+a_{DG}(\xi,\xi_{t}) & =(\eta_{tt},\xi_{t}),\end{align*}
 which implies \[
\cfrac{1}{2}\frac{d}{dt}\left(\|\xi_{t}\|_{L^{2}(\Omega)}^{2}+a_{DG}(\xi,\xi)\right)\leq\|\eta_{tt}\|_{L^{2}(\Omega)}\|\xi_{t}\|_{L^{2}(\Omega)}.\]
 Integrating from $t=0$ to $t=\tau$, we have \[
\begin{split}\|\xi_{t}(\tau)\|_{L^{2}(\Omega)}^{2}+ \frac{1}{2} \|\xi(\tau)\|_{a}^{2} & \leq\|\xi_{t}(0)\|_{L^{2}(\Omega)}^{2}+2\|\xi(0)\|_{a}^{2}+2\int_{0}^{\tau}\|\eta_{tt}\|_{L^{2}(\Omega)}\|\xi_{t}\|_{L^{2}(\Omega)}\\
 & \leq\|\xi_{t}(0)\|_{L^{2}(\Omega)}^{2}+2\|\xi(0)\|_{a}^{2}+2\max_{0\leq t\leq T}\|\xi_{t}\|_{L^{2}(\Omega)}\int_{0}^{T}\|\eta_{tt}\|_{L^{2}(\Omega)}.\end{split}
\]
 Therefore, we obtain \[
\|\xi_{t}\|_{L^{\infty}([0,T];L^{2}(\Omega))}^{2}+\|\xi\|_{L^{\infty}([0,T];a)}^{2}\leq C\Big(\|\xi_{t}(0)\|_{L^{2}(\Omega)}^{2}+\|\xi(0)\|_{a}^{2}+(\int_{0}^{T}\|\eta_{tt}\|_{L^{2}(\Omega)}\; dt)^{2}\Big).\]
 Finally, (\ref{eq:errorineq1}) is proved by noting that $\varepsilon=\eta-\xi$.

\begin{flushright}
$\square$
\par\end{flushright}

\begin{theorem} Let $\varepsilon,\eta$ and $\xi$ be the error quantities
defined in (\ref{eq:errordef}). Then we have the following error
bound \begin{equation}
\|\varepsilon\|_{L^{\infty}([0,T];L^{2}(\Omega))}\leq C\Big(\|\eta_{t}\|_{L^{1}([0,T],L^{2}(\Omega))}+\|\eta\|_{L^{\infty}([0,T];L^{2}(\Omega))}+ \mathcal{I}_2 \Big),
\label{eq:errorineq2}\end{equation}
where $\mathcal{I}_2 = \|\xi(0)\|_{L^{2}(\Omega)}$.
 \label{thm:est2} \end{theorem}

\textit{Proof:} Integrating by parts with respect to time in (\ref{eq:1}),
we have \[
-(\xi_{t},v_{t})+\partial_{t}(\xi_{t},v)+a_{DG}(\xi,v)=\partial_{t}(\eta_{t},v)-(\eta_{t},v_{t}).\]
 Taking $v(x,t)=\int_{t}^{\gamma}\xi(x,\tau)d\tau$, we have $v_{t}=-\xi$
and $v(\gamma)=0$. So, \[
(\xi_{t},\xi)-\partial_{t}(\xi_{t},v)-a_{DG}(v_{t},v)=\partial_{t}(\eta_{t},v)+(\eta_{t},\xi),\]
 which implies that \[
\frac{1}{2}\frac{d}{dt}\|\xi\|_{L^{2}(\Omega)}^{2}-\partial_{t}(\xi_{t},v)-\frac{1}{2}\frac{d}{dt}a_{DG}(v,v)=\partial_{t}(\eta_{t},v)+(\eta_{t},\xi).\]
 Integrating from $t=0$ to $t=\gamma$, we have \[
\cfrac{1}{2}\,\|\xi(\gamma)\|_{L^{2}(\Omega)}^{2}-\cfrac{1}{2}\,\|\xi(0)\|_{L^{2}(\Omega)}^{2}+(\xi_{t}(0),v(0))+\cfrac{1}{2}\, a_{DG}(v(0),v(0))=(\eta_{t}(0),v(0))+\int_{0}^{\gamma}(\eta_{t},\xi).\]
 Since $\xi_{t}-\eta_{t}=(u_{H}-u_h)_{t}$, we obtain \[
(\xi_{t}(0)-\eta_{t}(0),v(0))=((u_H-u_{h})_{t}(0),v(0))=0.\]
 Using the coercivity of $a_{DG}$, we have \begin{align*}
\|\xi(\gamma)\|_{L^{2}(\Omega)}^{2} & \leq\|\xi(0)\|_{L^{2}(\Omega)}^{2}+2\int_{0}^{\gamma}\|\eta_{t}\|_{L^{2}(\Omega)}\,\|\xi\|_{L^{2}(\Omega)}\\
 & \leq\|\xi(0)\|_{L^{2}(\Omega)}^{2}+2\max_{0\leq t\leq T}\|\xi\|_{L^{2}(\Omega)}\int_{0}^{T}\|\eta_{t}\|_{L^{2}(\Omega)}.\end{align*}
 Hence (\ref{eq:errorineq2}) is proved by noting that $\varepsilon=\eta-\xi$.

\begin{flushright}
$\square$
\par\end{flushright}

From Theorem \ref{thm:est1} and Theorem \ref{thm:est2}, we see that,
in order to estimate the error $\varepsilon=u_{h}-u_{H}$, we will need to find
a bound for $\eta$ given that the initial values $\xi_{t}(0)$ and
$\xi(0)$ are sufficiently accurate.

\subsection{Convergence analysis}

In this section, we will derive an error bound for $\eta=u_{h}-w_{H}$.
Notice that, on each coarse grid block $K$, we can express $u_{h}$
as \[
u_{h}=\sum_{i=1}^{n}c_{i,K}\widetilde{w}_{i,K}+\sum_{i=1}^{{n_0}}d_{i,K}z_{i,K}=u_{1,K}+u_{2,K}\]
 for some suitable coefficients $c_{i,K}$ and $d_{i,K}$ determined by a $L^2$-type projection,
 where $n_0$ is the dimension of $V_h^0(K)$.
 We write $u_h = u_1 + u_2$ with $u_i |_K = u_{i,K}$ for $i=1,2$.
 Moreover,
we recall that $C(u,f)$, defined in (\ref{eq:const}), is the constant
appearing in the consistency error estimate in Lemma \ref{lem:consistent}.
In the following theorem, we will give an estimate for the difference
between the fine grid solution $u_{h}$ and the projection of
$u_{h}$ into the coarse space $V_{H}$ defined in (\ref{eq:ellipticproj}).
The theorem says that such difference is bounded by a best approximation
error $\|u_{h}-v\|_{a}$ and a consistency error $hC(u,f)$. We emphasize
that, even though the coarse mesh size $H$ is fixed, but the fine
mesh size $h$ can be arbitrary small, and hence the consistency error
is small compared with the best approximation error $\|u_{h}-v\|_{a}$.

\begin{theorem} Let $w_{H}\in V_{H}$ be the solution of (\ref{eq:ellipticproj})
and $u_h$ be the solution of (\ref{eq:globalFEM}).
Then we have \begin{equation}
\|u_{h}-w_{H}\|_{a}\leq C(\|u_{h}-v\|_{a}+hC(u,f)),\quad\quad\forall\, v\in V_{H}.\label{eq:interbound}\end{equation}
 \label{thm:w-phi} \end{theorem}

%
%\marginpar{We need to say something about $C(u,f)$%
%}

\textit{Proof:} By the definition of $w_{H}$, we have \[
a_{DG}(w_{H},v)=a_{DG}(u_{h},v)+R_{u_{h}}(v),\quad\forall\, v\in V_{H}.\]
 So, we have \[
a_{DG}(w_{H}-v,w_{H}-v)=a_{DG}(u_{h}-v,w_{H}-v)+R_{u_{h}}(w_{H}-v).\]
 By (\ref{eq:coercivity}), (\ref{eq:continuity (bilinear)}) and
(\ref{eq:consistency}), we get \begin{align*}
\|w_{H}-v\|_{a}^{2} & \leq2a_{DG}(w_{H}-v,w_{H}-v)\\
 & =2a_{DG}(u_{h}-v,w_{H}-v)+2R_{u_{h}}(w_{H}-v)\\
 & \leq C(\|u_{h}-v\|_{a}+hC(u,f))\|w_{H}-v\|_{a}.\end{align*}
 Finally, we obtain \begin{align*}
\|u_{h}-w_{H}\|_{a} & \leq\|u_{h}-v\|_{a}+\|w_{H}-v\|_{a}\\
 & \leq C(\|u_{h}-v\|_{a}+hC(u,f)).
 \end{align*}

\begin{flushright}
$\square$
\par\end{flushright}

From the above theorem, we see that the error $\|u_{h}-w_{H}\|_{a}$
is controlled by the quantity $\|u_{h}-v\|_{a}$ for an arbitrary choice
of the function $v\in V_{H}$. Thus, to obtain our final error bound,
we only need to find a suitable function $v\in V_{H}$ to approximate
the finite element solution $u_{h}$. In the following theorem, we
will choose a specific $v$ in Theorem \ref{thm:w-phi} and prove
the corresponding error estimate.

\begin{theorem} Let $u_{h}\in V_{h}$ be the finite element solution.
Then we have
\begin{equation}
\|u_{h}-\phi\|^2_{a}\leq\sum_{K\in\mathcal{T}^{H}}
\Big (\frac{H}{\mu_{p+1,K}}(1+\cfrac{2a_{1}\gamma H}{h\mu_{p+1,K}})\int_{\partial K}(a\frac{\partial u_{1}}{\partial n})^{2}
+\cfrac{H^{2}}{\lambda_{m+1,K}}\|f-u_{tt}\|_{L^{2}(K)}^{2} \Big),
\label{eq:spectralconv}
\end{equation}
 where the function $\phi\in V_{H}$ is defined as \[
\phi|_{K}=\sum_{i=1}^{p}c_{i,K}\widetilde{w}_{i,K}+\sum_{i=1}^{m}d_{i,K}z_{i,K} = \phi_{1,K} + \phi_{2,K}.
\]
\label{thm:bound}
 \end{theorem}

\textit{Proof:}
For a given coarse grid block $K$, using the orthogonality condition (\ref{eq:Orthogonal}), we have \[
\int_{K}a|\nabla(u_{h}-\phi)|^{2}=\int_{K}a|\nabla(u_{1}-\phi_{1})|^{2}+\int_{K}a|\nabla(u_{2}-\phi_{2})|^{2}\]
 which implies \[
\|u_{h}-\phi\|_{a}^{2}=\|u_{1}-\phi_{1}\|_{a}^{2}+|u_{2}-\phi_{2}|_{a}^{2}, \]
where we write $\phi = \phi_1 + \phi_2$ and $\phi_i |_{K} = \phi_{i,K}$, for $i=1,2$.
 We will first estimate $\|u_{1}-\phi_{1}\|_{a}^{2}$. By the definition
of $a$-norm, we have \begin{equation}
\begin{split}\|u_{1}-\phi_{1}\|_{a}^{2} & = \sum_{K\in\mathcal{T}^{H}} \Big(\int_{K}a|\nabla(u_{1}-\phi_{1})|^{2}+\sum_{e\in\mathcal{E}^{H}}\cfrac{\gamma}{h}\int_{e}a|[(u_{1}-\phi_{1})]_{e}|^{2} \Big) \\
 & \leq\sum_{K\in\mathcal{T}^{H}}\left(\int_{K}a|\nabla(u_{1}-\phi_{1})|^{2}+\cfrac{2\gamma}{h}\int_{\partial K}a|(u_{1}-\phi_{1})|^{2}\right)\\
 & \leq\sum_{K\in\mathcal{T}^{H}}\left(\int_{K}a|\nabla(u_{1}-\phi_{1})|^{2}+\cfrac{2a_{1}\gamma}{h}\int_{\partial K}|(u_{1}-\phi_{1})|^{2}\right).\end{split}
\label{eq:first_est}\end{equation}
Next, we will estimate the right hand side of (\ref{eq:first_est}) for each $K$.

 We note that the eigenvalue problem (\ref{eq:eigenspace1}) is motivated
by the right hand side of (\ref{eq:first_est}). In particular, based
on the right hand side of (\ref{eq:first_est}), we consider \begin{equation}
\int_{K}a\nabla w_{\mu}\cdot\nabla v+\cfrac{1}{H}\int_{\partial K}w_{\mu}v=\widehat{\mu}\int_{K}R(w_{\mu})\cdot R(v),\quad\quad\forall\, v\in V_{H}^{1}(K),\label{eq:general_eig1}\end{equation}
 where the choice of $R$, e.g., $R=\sqrt{a}\nabla w_{\mu}$, depends
on how we would like to bound the error. Indeed, choosing the eigenvectors
that correspond to the largest $L_{K}$ eigenvalues, one can guarantee
that the best $L_{K}$ dimensional space \textit{in the space of snapshots}
is given by the first $L_{K}$ dominant eigenvectors. The choice of
$R(\cdot)$ is important and can influence the eigenvalue behavior.
For example, the use of oversampling domains both for the snapshot
space and the eigenvalue can provide a faster convergence. In this
paper, we take \[
R=\sqrt{a}\nabla w_{\mu},\]
 which allows estimating the right hand side of (\ref{eq:first_est})
by the energy norm. Note that, in (\ref{eq:eigenspace1}), we use
the smallest eigenvalues to determine the basis functions which is
the same as choosing the largest eigenvectors that correspond to the
largest eigenvalues of (\ref{eq:general_eig1}) because $\widehat{\mu}=1+\frac{1}{\mu}$.

Note that the eigenvalue problem (\ref{eq:eigenspace1}) is equivalent
to \[
a\frac{\partial w_{\mu}}{\partial n}=\frac{\mu}{H}w_{\mu}\quad\text{ on }\quad\partial K.\]
 So, for each $K$,
 \begin{equation}
\int_{\partial K}(a\frac{\partial u_{1,K}}{\partial n})^{2}
=\int_{\partial K}(a\frac{\partial}{\partial n}(\sum_{i=1}^{n}c_{i,K}w_{i,K}))^{2}=\int_{\partial K}(\sum_{i=1}^{n}\frac{\mu_{i,K}}{H}c_{i,K}w_{i,K})^{2}=\sum_{i=1}^{n}(\frac{\mu_{i,K}c_{i,K}}{H})^{2},
\label{eq:eigtmp}
\end{equation}
 where we have used the fact that $\int_{\partial K}w_{i,K}w_{j,K}=\delta_{ij}$.
Then,
by using the eigenvalue problem defined in (\ref{eq:eigenspace1}),
we have \[
\cfrac{1}{h}\int_{\partial K}|(u_{1,K}-\phi_{1,K})|^{2}=\cfrac{1}{h}\sum_{i=p+1}^{4n}c_{i,K}^{2}\leq\cfrac{H^{2}}{h\mu_{p+1,K}^{2}}\sum_{i=p+1}^{4n}\left(\frac{\mu_{i,K}}{H}\right)^{2}c_{i,K}^{2}\]
 and \[
\int_{K}a|\nabla(u_{1,K}-\phi_{1,K})|^{2}=\sum_{i=p+1}^{4n}\frac{\mu_{i,K}}{H}c_{i,K}^{2}\leq\frac{H}{\mu_{p+1,K}}\sum_{i=p+1}^{4n}(\frac{\mu_{i,K}}{H})^{2}c_{i,K}^{2}.\]
 Note that, by using (\ref{eq:eigtmp}), we have, \[
\sum_{i=p+1}^{4n}(\frac{\mu_{i,K}}{H})^{2}c_{i,K}^{2}\leq\sum_{i=1}^{4n}(\frac{\mu_{i,K}}{H})^{2}c_{i,K}^{2}=\int_{\partial K}(a\frac{\partial u_{1,K}}{\partial n})^{2}. \]
 Therefore \begin{equation}
\begin{split}\|u_{1}-\phi_{1}\|_{a}^{2}\leq &\sum_{K\in\mathcal{T}^{H}}
\Big( \frac{H}{\mu_{p+1,K}}(1+\cfrac{2a_{1}\gamma H}{h\mu_{p+1,K}})\sum_{i=p+1}^{4n}(\frac{\mu_{i,K}}{H})^{2}c_{i,K}^{2} \Big) \\
 \leq &\sum_{K\in\mathcal{T}^{H}} \Big( \frac{H}{\mu_{p+1,K}}(1+\cfrac{2a_{1}\gamma H}{h\mu_{p+1,K}})\int_{\partial K}(a\frac{\partial u_{1}}{\partial n})^{2}\Big). \end{split}
\end{equation}

Next, we will estimate $|u_{2}-\phi_{2}|_{a}^{2}$. Since $u_{h}$
satisfies \begin{align*}
\int_{K}a\nabla u_{h}\cdot\nabla v & =\int_{K}\left(f-(u_{h})_{tt}\right)v,\quad\quad\forall\, v\in V_h^0(K).\end{align*}
 Putting $v=z_{i,K}$, we obtain \[
\cfrac{\lambda_{i,K}}{H^{2}} \, d_{i,K}=\int_{K}a\nabla u_{h}\cdot \nabla z_{i,K}=\int_{K}\left(f-(u_{h})_{tt}\right) z_{i,K}.\]
 We define $f_{i,K}=\int_{K}\left(f-(u_{h})_{tt}\right) z_{i,K}$.
Then we have $f_{i,K}=\cfrac{\lambda_{i,K}}{H^{2}}d_{i,K}$ and \[
\sum_{i=1}^{n_0}f_{i,K}^{2}\leq\|f-(u_{h})_{tt}\|_{L^{2}(K)}^{2}.\]
 Hence, \begin{align*}
|u_{2}-\phi_{2}|_{a}^{2} & =\sum_{K\in\mathcal{T}^{H}}\int_{K}a|\nabla(u_{2}-\phi_{2})|^{2}\\
 & =\sum_{K\in\mathcal{T}^{H}}\sum_{i \geq m+1}\cfrac{\lambda_{i,K}}{H^{2}}d_{i,K}^{2}\\
 & \leq\sum_{K\in\mathcal{T}^{H}}\cfrac{H^{2}}{\lambda_{m+1,K}}\sum_{i \geq m+1}\cfrac{\lambda_{i,K}^{2}}{H^{4}}d_{i,K}^{2}\\
 & =\sum_{K\in\mathcal{T}^{H}}\cfrac{H^{2}}{\lambda_{m+1,K}}\sum_{i  \geq m+1}f_{i,K}^{2}\\
 & \leq\sum_{K\in\mathcal{T}^{H}}\cfrac{H^{2}}{\lambda_{m+1,K}}\|f-(u_{h})_{tt}\|_{L^{2}(K)}^{2}.\end{align*}

\begin{flushright}
$\square$
\par\end{flushright}

We note that, by the technique in \cite{eglms13}, we can also derive a bound for $\| u_1 - \phi_1 \|_a$ as follows
\begin{equation*}
\| u_1 - \phi_1 \|_a^2 \leq \sum_{K \in \mathcal{T}^H} \sum_{i \geq p+1} c_{i,K}^2.
\end{equation*}
This bound shows the decay of the error when more basis functions are used.

The bound in (\ref{eq:spectralconv}) gives the spectral convergence of
our GMsFEM.
Notice that, the term
\begin{equation*}
H \sum_{K\in\mathcal{T}^H} \sum_{\partial K} (a \frac{\partial u_1}{\partial n})^2
\end{equation*}
is uniformly bounded and can be considered as a norm for $u_1$.
Thus, (\ref{eq:spectralconv}) states that the error behaves as $O( \mu^{-1}_{p+1,K} + \lambda^{-1}_{m+1,K})$. We note that the eigenvalues increase
(and go to the infinity as the fine mesh size decreases) and thus
the error decreases as we increase the coarse space dimension.

Combining the results in Theorem \ref{thm:w-phi} and Theorem \ref{thm:bound},
we obtain
\begin{equation*}
\|\eta\|_{a}^{2}
\leq
C
\sum_{K\in\mathcal{T}^{H}}\left(\cfrac{H^{2}}{\lambda_{m+1,K}}\|f-u_{tt}\|_{L^{2}(K)}^{2}+\frac{H}{\mu_{p+1,K}}(1+\cfrac{2a_{1}\gamma H}{h\mu_{p+1,K}})\int_{\partial K}(a\frac{\partial u_{1}}{\partial n})^{2} \right)
+ h^2 C(u,f)^2.
\end{equation*}
Similarly, we obtain
\begin{equation*}
\|\eta_{t}\|_{a}^{2}
\leq
C \sum_{K\in\mathcal{T}^{H}}\left(\cfrac{H^{2}}{\lambda_{m+1,K}}\|f_{t}-u_{ttt}\|_{L^{2}(K)}^{2}+\frac{H}{\mu_{p+1,K}}(1+\cfrac{2a_{1}\gamma H}{h\mu_{p+1,K}})\int_{\partial K}(a\frac{\partial(u_{1})_{t}}{\partial n})^{2} \right)
+ h^2 C(u_t,f_t)^2.
\end{equation*}
Finally, using these bounds for $\eta$, as well as the estimates proved in Theorem
\ref{thm:est1} and Theorem \ref{thm:est2},
we obtain estimates for the error $\varepsilon$.

\section{Fully discretization }

\label{sec:full}

In this section, we will prove the convergence of the fully discrete
scheme (\ref{eq:fulldiscrete}). To simplify the notations, we define
the second order central difference operator $\delta^{2}$ by \[
\delta^{2}(u^{n})=\cfrac{u^{n+1}-2u^{n}+u^{n-1}}{\Delta t^{2}}.\]
 From the semi-discrete scheme (\ref{eq:ipdg}), we have \[
((u_{h})_{tt}^{n},v)+a_{DG}(u_{h}^{n},v)=(f^{n},v)-R_{u_{h}^{n}}(v)\]
 and the fully discrete scheme (\ref{eq:fulldiscrete}) can be written
as \[
(\delta^{2}(u_{H}^{n}),v)+a_{DG}(u_{H}^{n},v)=(f^{n},v),\quad\quad\mbox{for }n\geq1.\]
 Moreover, we define \begin{equation}
r^{n}=\begin{cases}
u_{tt}^{n}-\delta^{2}(w_{H}^{n}), & \mbox{for }n\geq1,\\
\Delta t^{-2}(\xi^{1}-\xi^{0}), & \mbox{for }n=0,\end{cases}\label{eq:res}\end{equation}
 and \[
R^{n}=\Delta t\sum_{i=0}^{n}r^{i}.\]

In order to prove the convergence for the fully discrete scheme, we
first prove the following lemma. The result will be needed in the
derivation of an upper bound for the time step size $\Delta t$.

%
%\marginpar{We can add 1-2 sentences about why we need this lemma%
%}
\begin{lemma} There exists a positive constant $\beta(h)$ such that
\[
a_{DG}(v,v)\leq \beta(h)^{-1}\|v\|_{L^{2}(\Omega)}^{2},\quad\forall\, v\in V_{H}.\]
Moreover, the constant $\beta(h)$ can be taken as $h^{2}a_{1}^{-1}(24+32\sqrt{3\Lambda}+16\gamma)^{-1}$.
\end{lemma}

\textit{Proof:} We first note that, if $p$ is a linear function defined
on the interval $I=[x_{1}-h/2,x_{1}+h/2]$, then we have \begin{align}
\|p\|_{L^{\infty}(I)}^{2} & \leq\cfrac{4}{h}\,\|p\|_{L^{2}(I)}^{2}\label{eq:infty to L2}\\
|p|_{H^{1}(I)}^{2} & \leq\cfrac{12}{h^{2}}\,\|p\|_{L^{2}(I)}^{2}.\label{eq:H1 to L2}\end{align}
 Then by the definition of $a_{DG}$ and the Cauchy-Schwarz inequality,
we have \begin{align*}
 & \: a_{DG}(v,v)\\
\leq & \:\sum_{K\in\mathcal{T}^{H}}\int_{K}a|\nabla v|^{2}-2\sum_{e\in\mathcal{E}^{H}}\Big(\int_{e}\{a\nabla v\cdot n\}_{e}\cdot[v]_{e}+\cfrac{\gamma}{h}\int_{e}a[v]_{e}^{2}\Big)\\
\leq & \: \sum_{K\in\mathcal{T}^{H}}\int_{K} a|\nabla v|^{2}+2(\sum_{K\in\mathcal{T}^{H}}h\|a\nabla v\cdot n_{\partial K}\|_{L^{2}(\partial K)}^{2})^{\frac{1}{2}}(\sum_{e\in\mathcal{E}^{H}}h^{-1}\|[v]\|_{L^{2}(e)}^{2})^{\frac{1}{2}}+\cfrac{\gamma}{h}\sum_{e\in\mathcal{E}^{H}}\int a[v]_{e}^{2}.\end{align*}
 Then by using (\ref{eq:discretetrace}), $a \leq a_1$ and estimating the jump terms
by $L^{2}(\partial K)$ norms, we have
\begin{equation}
\begin{split}
 & \: a_{DG}(v,v)\\
\leq & \: a_{1}\left(\sum_{K\in\mathcal{T}^{H}}\int_{K}|\nabla v|^{2}+2(\sum_{K\in\mathcal{T}^{H}} \Lambda\int_{K}|\nabla v|^{2})^{\frac{1}{2}}(\sum_{e\in\mathcal{E}^{H}}h^{-1}\|[v]\|_{L^{2}(e)}^{2})^{\frac{1}{2}}+\cfrac{\gamma}{h}\sum_{e\in\mathcal{E}^{H}}\int[v]_{e}^{2}\right)\\
\leq & \: a_{1}\left(\sum_{K\in\mathcal{T}^{H}}\int_{K}|\nabla v|^{2}+4(\sum_{K\in\mathcal{T}^{H}} \Lambda\int_{K}|\nabla v|^{2})^{\frac{1}{2}}(\sum_{K\in\mathcal{T}^{H}}h^{-1}\|v\|_{L^{2}(\partial K)}^{2})^{\frac{1}{2}}+\cfrac{2\gamma}{h}\sum_{K\in\mathcal{T}^{H}}\|v\|_{L^{2}(\partial K)}^{2}\right).
\end{split}
\label{eq:timestab}
\end{equation}
 Thus, it remains to estimate $\|\nabla v\|_{L^{2}(K)}$ and $\|v\|_{L^{2}(\partial K)}$.

We will estimate the term $\|\nabla v\|_{L^{2}(K)}$ first. For a
given coarse grid block $K$, we can write it as the union of fine grid
blocks $K=\cup_{F\subset K}F$, where we use $F$ to represent a generic
fine grid block. Since the fine grid blocks are rectangles, we can
write $F$ as a tensor product of two intervals, namely, $F=I_{x}^{F}\times I_{y}^{F}$.
For any $v\in V_{H}$ we can also write the restriction of $v$ on
$F$ as $v(x,y)=v_{F,1}(x)v_{F,2}(y)$. \begin{align*}
\int_{K}|\nabla v|^{2} & =\sum_{F\subset K}\int_{F}|\nabla v|^{2}\\
 & =\sum_{F\subset K}\left(h(v'_{F,2})^{2}\int_{I_{x}^{F}}(v_{F,1}(x))^{2}+h(v'_{F,1})^{2}\int_{I_{y}^{F}}(v_{F,2}(y))^{2}\right)\\
 & =\sum_{F\subset K}\left(\int_{I_{y}^{F}}(v'_{F,2}(y))^{2}\int_{I_{x}^{F}}(v_{F,1}(x))^{2}+\int_{I_{x}^{F}}(v'_{F,1}(x))^{2}\int_{I_{y}^{F}}(v_{F,2}(y))^{2}\right).\end{align*}
 Then, using (\ref{eq:H1 to L2}), we have \begin{align*}
\int_{K}|\nabla v|^{2} & \leq12h^{-2}\sum_{F\subset K}\left(\int_{I_{y}^{F}}(v{}_{F,2}(y))^{2}\int_{I_{x}^{F}}(v_{F,1}(x))^{2}+\int_{I_{x}^{F}}(v{}_{F,1}(x))^{2}\int_{I_{y}^{F}}(v_{F,2}(y))^{2}\right)\\
 & =24h^{-2}\sum_{F\subset K}\int_{F}|v|^{2}.\end{align*}

Next, we estimate the term $\|v\|_{L^{2}(\partial K)}$. For a generic
fine grid cell $F$, we write $I_{x}^{F}=[x_{1},x_{2}]$ and $I_{y}^{F}=[y_{1},y_{2}]$.
 Then, by using (\ref{eq:infty to L2}),
\begin{eqnarray*}
\|v\|_{L^{2}(\partial K)}^{2} & = & \sum_{F\subset K}\int_{\partial F\cap\partial K}(v_{F})^{2}\\
 & = & \sum_{F\subset K}\left(\int_{\partial F\cap(I_{x}\times\{y_{1}\})}(v{}_{F,2}(y_{1})v_{F,1}(x))^{2}+\int_{\partial F\cap(I_{x}\times\{y_{2}\})}(v{}_{F,2}(y_{2})v_{F,1}(x))^{2}\right)\\
 &  & +\sum_{F\subset K}\left(\int_{\partial F\cap(\{x_{1}\}\times I_{y})}(v{}_{F,1}(x_{1})v_{F,2}(y))^{2}+\int_{\partial F\cap(\{x_{2}\}\times I_{y})}(v{}_{F,1}(x_{2})v_{F,2}(y))^{2}\right)\\
 & \leq & \cfrac{4}{h}\sum_{F\subset K}\left(\int_{\partial F\cap(I_{x}\times\{y_{1}\})}\int_{[y_{1},y_{1}+h]}v{}_{F,2}(y)^{2}(v_{F,1}(x)^{2}+\int_{\partial F\cap(I_{x}\times\{y_{2}\})}\int_{[y_{2}-h,y_{2}]}v{}_{F,2}(y)^{2}v_{F,1}(x)^{2}\right)\\
 &  & +\cfrac{4}{h}\sum_{F\subset K}\left(\int_{\partial F\cap(\{x_{1}\}\times I_{y})}\int_{[x_{1},x_{1}+h]}v{}_{F,1}(x)^{2}v_{F,2}(y)^{2}+\int_{\partial F\cap(\{x_{2}\}\times I_{y})}\int_{[x_{2}-h,x_{2}]}v{}_{F,1}(x)v_{F,2}(y)^{2}\right)\\
 & \leq & \cfrac{4}{h}\sum_{F\subset K}\left(\int_{F}(v{}_{F,2}(y)v_{F,1}(x))^{2}+\int_{F}(v{}_{F,1}(x)v_{F,2}(y))^{2}\right)\\
 & = & \cfrac{8}{h} \, \|v\|_{L^{2}(K)}^{2}.\end{eqnarray*}

% Then, by using (\ref{eq:infty to L2}), we get \begin{eqnarray*}
%\|v\|_{L^{2}(\partial K)}^{2} & \leq & \cfrac{4}{h}\sum_{F\subset K}\left(\int_{\partial F\cap(I_{x}\times\{y_{1}))}\int_{[y_{1},y_{1}+h]}v{}_{F,2}(y)^{2}v_{F,1}(x)^{2}+\int_{\partial F\cap(I_{x}\times\{y_{2}))}\int_{[y_{2}-h,y_{2}]}v{}_{F,2}(y)^{2}v_{F,1}(x)^{2}\right)\\
 %&  & +\cfrac{4}{h}\sum_{F\subset K}\left(\int_{\partial F\cap(\{x_{1}\}\times I_{y})}\int_{[x_{1},x_{1}+h]}v{}_{F,1}(x)^{2}v_{F,2}(y)^{2}+\int_{\partial F\cap(\{x_{2}\}\times I_{y})}\int_{[x_{2}-h,x_{2}]}v{}_{F,1}(x)v_{F,2}(y)^{2}\right)\\
 %& \leq & \cfrac{4}{h}\sum_{F\subset K}\left(\int_{F}(v{}_{F,2}(y)v_{F,1}(x))^{2}+\int_{F}(v{}_{F,1}(x)v_{F,2}(y))^{2}\right)\\
 %& = & \cfrac{8}{h}\|v\|_{L^{2}(K)}^{2}.\end{eqnarray*}

Consequently, combining the above results and using (\ref{eq:timestab}),
\begin{equation*}
a_{DG}(v,v) \leq
\frac{a_1}{h^2} \Big( 24 + 32 \sqrt{3\Lambda} + 16\gamma \Big) \|v\|_{L^2(\Omega)}^2.
\end{equation*}

%\begin{align*}
%a_{DG}(v,v) & \leq a_{1}\left(\sum_{K\in\mathcal{T}^{H}}\int_{K}|\nabla v|^{2}+2(\sum_{K\in\mathcal{T}^{H}}\int_{K}|\nabla v|^{2})^{\frac{1}{2}}(\sum_{e\in\mathcal{E}^{H}}\cfrac{\|[v]\|_{L^{2}(e)}^{2}}{h})^{\frac{1}{2}}+\cfrac{\gamma}{h}\sum_{e\in\mathcal{E}^{H}}\int[v]_{e}^{2}\right)\\
 %& \leq\cfrac{a_{1}}{h^{2}}(24+32\sqrt{3}+16\gamma)\|v\|_{L^{2}(\Omega)}^{2}.\end{align*}

\begin{flushright}
$\square$
\par\end{flushright}

Finally, we will state and prove
the convergence of the fully discrete scheme (\ref{eq:fulldiscrete}).

\begin{theorem} \label{thm:timestep} Assume that the time step size $\Delta t$
satisfies the stability condition $\Delta t^{2}<4\beta(h)$. We have \begin{equation}
\max_{0\leq n\leq N}\|\varepsilon^{n}\|_{L^{2}(\Omega)}\leq C\left(\|\varepsilon^{0}\|_{L^{2}(\Omega)}+\max_{0\leq n\leq N}\|\eta^{n}\|_{L^{2}(\Omega)}+\Delta t\sum_{n=0}^{N}\|R^{n}\|_{L^{2}(\Omega)}\right).\label{eq:dis L2 error}\end{equation}
 \end{theorem}

\textit{Proof:} Notice that \[
(\delta^{2}(u_{H}^{n}-w_{H}^{n}+w_{H}^{n}-u_{h}^{n}),v)+a_{DG}(u_{H}^{n}-u_{h}^{n},v)=((u_{h})_{tt}^{n}-\delta^{2}(u_{h}^{n}),v)+R_{u_{h}^{n}}(v).\]
 Using the definitions of $\xi$ and $w_{H}$, we have \[
(\delta^{2}(\xi^{n}),v)+a_{DG}(\xi^{n},v)=(r^{n},v),\quad\quad\mbox{for }n\geq1.\]
 So we have \[
(\cfrac{\xi^{n+1}-\xi^{n}}{\Delta t},v)-(\cfrac{\xi^{n}-\xi^{n-1}}{\Delta t},v)+\Delta t\, a_{DG}(\xi^{n},v)=\Delta t(r^{n},v).\]
 Summing up, for $n\geq1$, \[
(\cfrac{\xi^{n+1}-\xi^{n}}{\Delta t},v)-(\cfrac{\xi^{1}-\xi^{0}}{\Delta t},v)+\Delta t\sum_{i=1}^{n}a_{DG}(\xi^{i},v)=\Delta t\sum_{i=1}^{n}(r^{i},v).\]
 To simplify the notations, we define \[
\Xi^{n}=\Delta t\sum_{i=1}^{n}\xi^{i},\quad\text{ for }n\geq1;\quad\text{and}\quad\Xi^{0}=0.\]
 Then we get \[
\left(\cfrac{\xi^{n+1}-\xi^{n}}{\Delta t},v\right)+a_{DG}(\Xi^{n},v)=(R^{n},v),\quad n\geq1.\]
 Substituting $v=\xi^{n+1}+\xi^{n}$, we have \[
\|\xi^{n+1}\|^{2}_{L^2(\Omega)} - \|\xi^{n}\|^{2}_{L^2(\Omega)} +\Delta t\, a_{DG}(\Xi^{n},\xi^{n+1}+\xi^{n})=\Delta t(R^{n},\xi^{n+1}+\xi^{n}),\]
 and summing for all $n\geq1$, we have \[
\|\xi^{n+1}\|^{2}_{L^2(\Omega)} - \|\xi^{1}\|^{2}_{L^2(\Omega)} +\Delta t\sum_{i=1}^{n}a_{DG}(\Xi^{i},\xi^{i+1}+\xi^{i})=\Delta t\sum_{i=1}^{n}(R^{i},\xi^{i+1}+\xi^{i}).\]
 Notice that we have $\Xi^{n+1}-\Xi^{n-1}=\Delta t(\xi^{n+1}+\xi^{n})$
for $n\geq1$. So \begin{align*}
\Delta t\sum_{i=1}^{n}\tilde{a}_{DG}(\Xi^{i},\xi^{i+1}+\xi^{i}) & =\sum_{i=1}^{n}a_{DG}(\Xi^{i},\Xi^{i+1}-\Xi^{i-1})\\
 & =\sum_{i=1}^{n}a_{DG}(\Xi^{i},\Xi^{i+1})-\sum_{i=0}^{n-1}a_{DG}(\Xi^{i},\Xi^{i+1})\\
 & =a_{DG}(\Xi^{n},\Xi^{n+1}).\end{align*}
 Moreover, \begin{align*}
a_{DG}(\Xi^{n},\Xi^{n+1}) & =a_{DG}(\cfrac{\Xi^{n}+\Xi^{n+1}}{2},\cfrac{\Xi^{n}+\Xi^{n+1}}{2})-a_{DG}(\cfrac{\Xi^{n}-\Xi^{n+1}}{2},\cfrac{\Xi^{n}-\Xi^{n+1}}{2})\\
 & \geq-\cfrac{\Delta t^{2}}{4}a_{DG}(\xi^{n+1},\xi^{n+1}).\end{align*}
 So we have \[
\|\xi^{n+1}\|^{2}_{L^2(\Omega)} -\cfrac{\Delta t^{2}}{4}a_{DG}(\xi^{n+1},\xi^{n+1})\leq\|\xi^{1}\|^{2}_{L^2(\Omega)}+\Delta t\sum_{i=1}^{n}(R^{i},\xi^{i+1}+\xi^{i}),\quad n\geq1.\]

Using the assumption $\Delta t^{2}<4\beta(h)$, we have $C_{s}=1-\cfrac{\Delta t^{2}}{4\beta(h)}>0$.
Therefore, \begin{align*}
C_{s}\|\xi^{n+1}\|^{2}_{L^2(\Omega)} & \leq\|\xi^{1}\|^{2}_{L^2(\Omega)}+\Delta t\sum_{i=1}^{n}(R^{i},\xi^{i+1}+\xi^{i})\\
 & \leq\|\xi^{1}\|^{2}_{L^2(\Omega)}+2\Delta t\max_{1\leq i\leq n+1}\{\|\xi^{i}\|_{L^{2}(\Omega}\}\sum_{i=1}^{n}\|R^{i}\|_{L^{2}(\Omega)}\\
 & \leq\|\xi^{1}\|^{2}_{L^2(\Omega)}+\cfrac{C_{s}}{2}\max_{1\leq i\leq n+1}\{\|\xi^{i}\|_{L^{2}(\Omega)}\}^{2}+\cfrac{2}{C_{s}}\left(\Delta t\sum_{i=1}^{n}\|R^{i}\|_{L^{2}(\Omega)}\right)^{2}.\end{align*}
 Then \[
\max_{1\leq i\leq n+1}\{\|\xi^{i}\|_{L^{2}(\Omega)}\}\leq\sqrt{\cfrac{2}{C_{s}}}\|\xi^{1}\|_{L^{2}(\Omega)}+\cfrac{2}{C_{s}}\Delta t\sum_{i=1}^{n}\|R^{i}\|_{L^{2}(\Omega)}.\]
 Since $\xi^{1}=\xi^{0}+\Delta t^{2}r^{0}$, we have \[
\max_{1\leq i\leq n+1}\{\|\xi^{i}\|_{L^{2}(\Omega)}\}\leq C\left(\|\xi^{0}\|_{L^{2}(\Omega)}+\Delta t^{2}\|r^{0}\|_{L^{2}(\Omega)}+\Delta t\sum_{i=1}^{n}\|R^{i}\|_{L^{2}(\Omega)}\right)\]
 and using the definition of $R^{0}$, \[
\max_{1\leq i\leq n+1}\{\|\xi^{i}\|_{L^{2}(\Omega)}\}\leq C\left(\|\xi^{0}\|_{L^{2}(\Omega)}+\Delta t\sum_{i=0}^{n}\|R^{i}\|_{L^{2}(\Omega)}\right).\]
 Thus, \[
\max_{0\leq i\leq n+1}\{\|\xi^{i}\|_{L^{2}(\Omega)}\}\leq C\left(\|\xi^{0}\|_{L^{2}(\Omega)}+\Delta t\sum_{i=0}^{n}\|R^{i}\|_{L^{2}(\Omega)}\right).\]
Finally, by using the relation $\varepsilon = \eta - \xi$, we obtain (\ref{eq:dis L2 error}).

\begin{flushright}
$\square$
\par\end{flushright}

Next, we will estimate the right hand side of (\ref{eq:dis L2 error}).
To do so, we prove the following lemmas.

\begin{lemma} We have \[
\|r^{0}\|_{L^{2}(\Omega)}\leq C(\Delta t^{-1}\|\eta_{t}\|_{L^{\infty}([0,T];L^{2}(\Omega))}+\Delta t\|(u_{h})_{ttt}\|_{C([0,T];L^{2}(\Omega))}).\]
 \end{lemma}

\textit{Proof:} By (\ref{eq:res}), we have $r^{0}=\Delta t^{-2}(\xi^{1}-\xi^{0})$
and by the definition of $u_{H}^{0}$, we have \begin{align*}
(u_{H}^{0}-u^{0},v) & =0,\quad\forall\, v\in V_{H}.\end{align*}
 Then using the definitions of $\xi^{1}$ and $\xi^{0}$, we have
\begin{align*}
(\xi^{1}-\xi^{0},v) & =(u_{H}^{1}-w_{H}^{1},v)-(u_{H}^{0}-w_{H}^{0},v)\\
 & =(u_{h}^{1}-w_{H}^{1},v)+(u_{H}^{1}-u_{h}^{1},v)-(u_{h}^{0}-w_{H}^{0},v)\\
 & =((u_{h}^{1}-u_{h}^{0})-(w_{H}^{1}-w_{H}^{0}),v)+(u_{H}^{1}-u_{h}^{1},v).\end{align*}
 The first term can be estimated in the following way \begin{align*}
|((u_{h}^{1}-u_{h}^{0})-(w_{H}^{1}-w_{H}^{0}),v)| & \leq\Big|(\int_{0}^{t_{1}}\partial_{t}(u_{h}-w_{H}),v)\Big|\\
 & \leq\Delta t\,\|\eta_{t}\|_{L^{\infty}([0,T];L^{2}(\Omega))}\|v\|_{L^{2}(\Omega)}.\end{align*}
 To estimate the second term, by the Taylor's expansion, we get \[
u_{h}^{1}=u_{h}^{0}+\Delta t\,(u_{h})_{t}^{0}+\cfrac{\Delta t^{2}}{2}(u_{h})_{tt}^{0}+\cfrac{\Delta t^{3}}{6}(u_{h})_{ttt}(\cdot,s),\quad\quad\mbox{where }0<s<t^{1}.\]
 By the definition of $u_{H}^{1}$, \[
(u_{H}^{1},v)=(u_{h}^{1},v)=(u_{h}^{0}+\Delta t(u_{h})_{t}^{0}+\cfrac{\Delta t^{2}}{2}(\tilde{v},v)\]
 Thus, \begin{align*}
(u_{H}^{1}-u_{h}^{1},v) & =\cfrac{\Delta t^{2}}{2}(\tilde{v}-(u_{h})_{tt}^{0},v)-\cfrac{\Delta t^{3}}{6}((u_{h})_{ttt}(\cdot,s),v)\\
 & =\cfrac{\Delta t^{2}}{2}[(f^{0},v)-a(u_{h}^{0},v)+((u_{h})_{tt}^{0},v)]-\cfrac{\Delta t^{3}}{6}((u_{h})_{ttt}(\cdot,s),v)\\
 & =-\cfrac{\Delta t^{3}}{6}((u_{h})_{ttt}(\cdot,s),v)\end{align*}
 which proves the Lemma.

\begin{flushright}
$\square$
\par\end{flushright}

\begin{lemma} For $n\geq1$, we have \[
\|r^{n}\|_{L^{2}(\Omega)}\leq C(\Delta t^{-1}\int_{t_{n-1}}^{t_{n+1}}\|\eta_{tt}(\cdot,\tau)\|_{L^{2}(\Omega)}+\Delta t\int_{t_{n-1}}^{t_{n+1}}\|(u_{h})_{tttt}(\cdot,\tau)\|_{L^{2}(\Omega)}).\]
 \end{lemma}

\textit{Proof:} By the definition of $r^{n}$, \begin{align*}
\|r^{n}\|_{L^{2}(\Omega)} & =\|(u_{h})_{tt}^{n}-\delta^{2}w_{H}^{n}\|_{L^{2}(\Omega)}\\
 & \leq\|\delta^{2}(w_{H}^{n}-u_{h}^{n})\|_{L^{2}(\Omega)}+\|(u_{h})_{tt}^{n}-\delta^{2}u_{h}^{n}\|_{L^{2}(\Omega)}.\end{align*}
 Using the identity \[
v^{n+1}-2v^{n}+v^{n-1}=\Delta t\int_{t_{n-1}}^{t_{n+1}}\left(1-\cfrac{|\tau-t_{n}|}{\Delta t}\right)v_{tt}(\tau)d\tau,\]
 the first term can be estimated as follows \begin{align*}
(\delta^{2}(w_{H}^{n}-u_{h}^{n}),v) & =\cfrac{1}{\Delta t}\int_{t_{n-1}}^{t_{n+1}}(1-\cfrac{|\tau-t_{n}|}{\Delta t})\left(((w_{H})_{tt}-(u_{h})_{tt}),v\right)(\tau)d\tau\\
 & \leq\cfrac{1}{\Delta t}\int_{t_{n-1}}^{t_{n+1}}\|\eta_{tt}(\cdot,\tau)\|_{\tilde{L}^{2}(\Omega)}\,\|v\|_{\tilde{L}^{2}(\Omega)}d\tau.\end{align*}
 To estimate the term $\|(u_{h})_{tt}^{n}-\delta^{2}u_{h}^{n}\|_{L^{2}(\Omega)}$,
we use \[
\delta^{2}u_{h}^{n}=(u_{h})_{tt}^{n}+\cfrac{1}{6\Delta t^{2}}\int_{t_{n-1}}^{t_{n+1}}(\Delta t-|\tau-t_{n}|)^{3}(u_{h})_{tttt}(\cdot,\tau)d\tau.\]
 This implies \[
\|(u_{h})_{tt}^{n}-\delta^{2}u_{h}^{n}\|_{L^{2}(\Omega)}\leq\cfrac{\Delta t}{6}\int_{t_{n-1}}^{t_{n+1}}\|(u_{h})_{tttt}(\cdot,\tau)\|_{L^{2}(\Omega)}d\tau.\]

\begin{flushright}
$\square$
\par\end{flushright}

Using the definition of $R^n$ and the above two lemma, we get
\begin{equation*}
\begin{split}
&\: \|R^{n}\|_{L^{2}(\Omega)} \\
\leq &\: C\Big(\int_{0}^{t_{n}}\|\eta_{tt}(\cdot,\tau)\|_{L^2(\Omega)}+\|\eta_{t}\|_{L^\infty([0,T];L^2(\Omega))}+\Delta t^{2}\int_{0}^{t_{n}}\|(u_{h})_{tttt}(\cdot,\tau)\|_{L^{2}(\Omega)}+\Delta t^{2}\|(u_{h})_{ttt}\|_{C([0,T];L^{2}(\Omega))} \Big).
\end{split}
\end{equation*}
 Hence we obtain
 \begin{align*}
&\: \Delta t\sum_{n=0}^N \|R^{n}\|_{L^{2}(\Omega)} \\
 \leq &\: 2T\max_{0\leq n \leq N}\|R^{n}\|_{L^{2}(\Omega)}\\
  \leq &\: C(\int_{0}^{T}\|\eta_{tt}(\cdot,\tau)\|_{L^2(\Omega)}+\|\eta_{t}\|_{L^\infty([0,T];L^2(\Omega))}+\Delta t^{2}\int_{0}^{T}\|(u_{h})_{tttt}(\cdot,\tau)\|_{L^{2}(\Omega)}+\Delta t^{2}\|(u_{h})_{ttt}\|_{C([0,T];L^{2}(\Omega))}). \end{align*}
 Combining the estimates of $\eta$ proved in Section 3 and (\ref{eq:dis L2 error}),
we obtain the error estimate for the fully discrete scheme (\ref{eq:fulldiscrete}).

\section{Conclusions}

In this paper, we present a multiscale simulation method based on
Generalized Multiscale Finite Element Method for solving the wave
equation in heterogeneous media. For the construction of multiscale
basis functions, we divide the snapshot space into two spaces.
The first snapshot space represents the degrees of the freedom
associated with boundary nodes and consists of $a$-harmonic functions.
The second snapshot space represents the interior degrees of the freedom
and consists of all zero Dirichlet vectors. For each  snapshot
space, we introduce local spectral problems motivated by the analysis
presented in the paper.
We use these local spectral problems to identify important
modes in each of the snapshot spaces.
The local spectral problems are designed to
achieve a high accuracy and motivated by the global coupling formulation.
The use of
multiple snapshot spaces
and multiple spectral problems is one of the novelties of this work.
Using the dominant modes from local spectral problems,
multiscale basis functions are constructed to
represent the solution space locally within each coarse block.
These multiscale basis functions
are coupled via the symmetric interior penalty discontinuous Galerkin
method which provides a block diagonal mass matrix, and, consequently,
results in fast computations in an explicit time discretization.
Numerical examples are presented.  In particular,
we discuss how the modes from our snapshot spaces can affect
the accuracy
of the method. Our numerical results show that
one can obtain an accurate approximation of the solution with  GMsFEM
using less than
$3\%$ of the total local degrees of freedom.
We also test oversampling strategies following \cite{eglp13}.
 Analysis of the method is presented.

\section{Appendix}

In this Appendix, we will prove Lemma \ref{lem:consistent}. Let $v\in V_{H}^{1}$.
By assumption, $u\in H^{2}(\Omega)$, thus $a_{DG}(u,v)$ is well-defined
and we have \[
\sum_{K\in\mathcal{T}^{H}}(\frac{\partial^{2}u}{\partial t^{2}},v)_{L^{2}(K)}+a_{DG}(u,v)=\sum_{K\in\mathcal{T}^{H}}(f,v)_{L^{2}(K)},\quad\forall\, v\in H^{1}(\mathcal{T}^{H}).\]
 Moreover, the following standard finite element error estimate holds
\[
|u-u_{h}|_{H^{1}(\Omega)}\leq Ch|u|_{H^{2}(\Omega)}.\]
 By the definition of the consistency error, we have \begin{equation}
R_{u_{h}}(v)=(\frac{\partial^{2}u}{\partial t^{2}}-\frac{\partial^{2}u_{h}}{\partial t^{2}},v)+a(u,v)-a_{DG}(u_{h},v),\quad\quad\forall\, v\in H^{1}(\mathcal{T}^{H}).\label{eq:consistent 1}\end{equation}
 %And by (\ref{eq:globalFEM}) \begin{align}
%(f,v)-(\frac{\partial^{2}u_{FEM}}{\partial t^{2}},v)-a_{DG}(u_{FEM},v) & =0\;\forall v\in V_{h}\label{eq:consistent 2}\end{align}
%For $v\in V^{H}$, we define $v_{c}\in V_{h}$ such that \begin{align*}
%v_{c}(x) & =\cfrac{\left(\sum_{K\ni x}v_{K}(x)\right)}{\left(\sum_{K\ni x}1\right)}\mbox{ for all nodal point }x.\end{align*}
Next, we define $v_{c}\in V_{h}$ in the following way. For each vertex
in the triangulation, the value of $v_{c}$ is defined as the average
value of $v$ at this vertex. Then by direct calculations, we have
\[
\sum_{K\in\mathcal{T}^{H}}|v-v_{c}|_{H^{1}(K)}^{2}\leq C\cfrac{1}{h}\sum_{e\in\mathcal{E}^{H}}\|[v]\|_{L^{2}(e)}^{2}\]
 and \[
\sum_{K\in\mathcal{T}^{H}}\|v-v_{c}\|_{L^{2}(K)}^{2}\leq Ch\sum_{e\in\mathcal{E}^{H}}\|[v]\|_{L^{2}(e)}^{2}.\]
 Clearly, we have $[v-v_{c}]_{e}=[v]_{e}$ for all $e\in\mathcal{E}^{H}$
since $v_{c}\in C^{0}(\Omega)$. Therefore we get \[
\|v-v_{c}\|_{H^{1}(\mathcal{T}^{H})}^{2}\leq C\cfrac{1}{h}\sum_{e\in\mathcal{E}^{H}}\|[v]\|_{L^{2}(e)}^{2}.\]

By (\ref{eq:consistent 1}) and (\ref{eq:globalFEM}) as well as the
fact that $a_{DG}(u_{h},v_{c})=a(u_{h},v_{c})$, we have \begin{align*}
R_{u_{h}}(v) & =\sum_{K\in\mathcal{T}^{H}}(\frac{\partial^{2}(u-u_{h})}{\partial t^{2}},v-v_{c})_{L^{2}(K)}+a_{DG}(u-u_{h},v-v_{c}).\end{align*}
 Next, we will estimate the two terms on the right hand side. For the
first term, we have \begin{align*}
\sum_{K\in\mathcal{T}^{H}}(\frac{\partial(u-u_{h})}{\partial t^{2}},v-v_{c})_{L^{2}(K)} & \leq\|\frac{\partial^{2}(u-u_{h})}{\partial t^{2}}\|_{L^{2}(\Omega)}\left(\sum_{K\in\mathcal{T}^{H}}\|v-v_{c}\|_{L^{2}(K)}^{2}\right)^{\frac{1}{2}}\\
 & \leq Ch\,\|\frac{\partial^{2}(u-u_{h})}{\partial t^{2}}\|_{L^{2}(\Omega)}\left(\cfrac{1}{h}\sum_{e\in\mathcal{E}^{H}}\|[v]\|_{L^{2}(e)}^{2}\right)^{\frac{1}{2}}.\end{align*}
 For the second term, by the definition of $a_{DG}$ and the Cauchy-Schwarz
inequality, we have \begin{align*}
 & \: a_{DG}(u-u_{h},v-v_{c})\\
= & \:\sum_{K\in\mathcal{T}^{H}}\int_{K}a\nabla(u-u_{h})\cdot\nabla(v-v_{c})-\sum_{e\in\mathcal{E}^{H}}\int_{e}\{a\nabla(u-u_{h})\cdot n\}[v]\\
\leq & \:\sum_{K\in\mathcal{T}^{H}}a_{1}|u-u_{h}|_{H^{1}(K)}\;|v-v_{c}|_{H^{1}(K)}+\left(\sum_{K\in\mathcal{T}^{H}}\int_{\partial K}(a\nabla(u-u_{h})\cdot n)^{2}\right)^{\frac{1}{2}}\left(\sum_{e\in\mathcal{E}^{H}}\int_{e}[v]^{2}\right)^{\frac{1}{2}}.\end{align*}
 To estimate the flux term above, we let $I_{K}$ be the standard
finite element interpolant. %Then we will estimate the term $\left(\sum_{K\in\mathcal{T}^{H}}\int_{\partial K}(a\nabla(u-u_{FEM})\cdot n)^{2}\right)$.
Then we have \begin{align*}
 & \:\sum_{K\in\mathcal{T}^{H}}\int_{\partial K}(a\nabla(u-u_{h})\cdot n)^{2}\\
\leq & \:2\left(\sum_{K\in\mathcal{T}^{H}}\int_{\partial K}(a\nabla(u-I_{K}(u))\cdot n)^{2}\right)+2\left(\sum_{K\in\mathcal{T}^{H}}\int_{\partial K}(a\nabla(I_{K}(u)-u_{h})\cdot n)^{2}\right)\\
\leq & \: Ca_{1}\left(h|u|_{H^{2}(K)}^{2}+\cfrac{1}{h}|I_{K}(u)-u_{h}|_{H^{1}(K)}^{2}\right)\\
\leq & \: Ca_{1}\left(h|u|_{H^{2}(K)}^{2}+\cfrac{1}{h}|I_{K}(u)-u|_{H^{1}(K)}^{2}+\cfrac{1}{h}|u-u_{h}|_{H^{1}(K)}^{2}\right)\\
\leq & \: Ca_{1}h|u|_{H^{2}(K)}^{2}.\end{align*}
 Next we will estimate the term $\sum_{K\in\mathcal{T}^{H}}a_{1}|u-u_{h}|_{H^{1}(K)}\;|v-v_{c}|_{H^{1}(K)}.$
We have \begin{align*}
 & \:\sum_{K\in\mathcal{T}^{H}}a_{1}|u-u_{h}|_{H^{1}(K)}\cdot|v-v_{c}|_{H^{1}(K)}\\
\leq & \: Ca_{1}(\sum_{K\in\mathcal{T}^{H}}|u-u_{h}|_{H^{1}(K)}^{2})^{\frac{1}{2}}(\sum_{K\in\mathcal{T}^{H}}|v-v_{c}|_{H^{1}(K)}^{2})^{\frac{1}{2}}\\
\leq & \: Ca_{1}(\cfrac{1}{h}\sum_{K\in\mathcal{T}^{H}}|u-u_{h}|_{H^{1}(K)}^{2})^{\frac{1}{2}}(\sum_{e\in\mathcal{E}^{H}}\|[v]\|_{L^{2}(e)}^{2})^{\frac{1}{2}}\\
\leq & \: Ca_{1}h|u|_{H^{2}}(\frac{1}{h}\sum_{e\in\mathcal{E}^{H}}\|[v]\|_{L^{2}(e)}^{2})^{\frac{1}{2}}.\end{align*}
 Combining the above estimates, we get \[
|R_{u_{h}}(v)|\leq\cfrac{C}{a_{0}}h\left(\|(u-u_{h})_{tt}\|_{L^{2}(K)}+a_{1}|u|_{H^{2}(\Omega)}\right)\left(\cfrac{1}{h}\sum_{e\in\mathcal{E}^{H}}\|a[v]\|_{L^{2}(e)}^{2}\right)^{\frac{1}{2}}.\]

Finally, we assume that the second time derivatives of $u$ and $u_{h}$
are smooth functions. Then we have \[
(u_{ttt},v)+\int_{\Omega}a\nabla u_{t}\cdot\nabla v=\int_{\Omega}f_{t}v.\]
 Letting $v=u_{tt},$ we have \[
\frac{d}{dt}(\|u_{tt}\|_{L^{2}(\Omega)}^{2}+\int_{\Omega}a|\nabla u_{t}|^{2})=\int_{\Omega}f_{t}u_{tt}\]
 which leads to \[
\|u_{tt}\|_{L^{\infty}([0,T];L^{2}(\Omega))}+\|u_{t}\|_{L^{\infty}([0,T];a)}\leq C(\|f_{t}\|_{L^{1}([0,T];L^{2}(\Omega))}+\|u_{tt}(\cdot,0)\|_{L^{2}(\Omega)}+|u_{t}(\cdot,0)|_{a}).\]
 Similarly, for the finite element solution $u_{h}$, we have \begin{align*}
 & \:\|(u_{h})_{tt}\|_{L^{\infty}([0,T];L^{2}(\Omega))}+\|(u_{h})_{t}\|_{L^{\infty}([0,T];a)}\\
\leq & \: C(\|f_{t}\|_{L^{1}([0,T];L^{2}(\Omega))}+\|(u_{h})_{tt}(\cdot,0)\|_{L^{2}(\Omega)}+|(u_{h})_{t}(\cdot,0)|_{a}).\end{align*}
 Consequently, we get \[
|R_{u_{h}}(v)|\leq Ch(\|f_{t}\|_{L^{1}([0,T];L^{2}(\Omega))}+\|u_{tt}(\cdot,0)\|_{L^{2}(\Omega)}+a_{1}|u|_{H^{2}(\Omega)})\|v\|_{H^{1}(\mathcal{T}^{H})}.\]
 Finally, the constant $C(u,f)$ in the lemma can be chosen as \begin{equation}
C(u,f)\approx\|f_{t}\|_{L^{1}([0,T];L^{2}(\Omega))}+\|u_{tt}(\cdot,0)\|_{L^{2}(\Omega)}+a_{1}|u|_{H^{2}(\Omega)}.\label{eq:const}\end{equation}

\bibliographystyle{plain}
\bibliography{segabs_ceg,EricPaper,EricRef}

\begin{thebibliography}{10}

\bibitem{dispersion}
H.~Chan, E.~Chung, and G.~Cohen.
\newblock Stability and dispersion analysis of staggered discontinuous
  {G}alerkin method for wave propagation.
\newblock {\em Int. J. Numer. Anal. Model.}, 10:233--256, 2013.

\bibitem{AADA}
E.~Chung, Y.~Efendiev, and R.~Gibson.
\newblock An energy-conserving discontinuous multiscale finite element method
  for the wave equation in heterogeneous media.
\newblock {\em Advances in Adaptive Data Analysis}, 3:251--268, 2011.

\bibitem{ChungEngquist06}
E.~Chung and B.~Engquist.
\newblock Optimal discontinuous {G}alerkin methods for wave propagation.
\newblock {\em SIAM J. Numer. Anal.}, 44:2131--2158, 2006.

\bibitem{ChungEngquist09}
E.~Chung and B.~Engquist.
\newblock Optimal discontinuous {G}alerkin methods for the acoustic wave
  equation in higher dimensions.
\newblock {\em SIAM J. Numer. Anal.}, 47:3820--3848, 2009.

\bibitem{JCAM-meta}
E.~Chung and P.~Ciarlet Jr.
\newblock A staggered discontinuous {G}alerkin method for wave propagation in
  media with dielectrics and meta-materials.
\newblock {\em J. Comput. Appl. Math.}, 239:189--207, 2013.

\bibitem{JCP-max}
E.~Chung, P.~Ciarlet Jr., and T.~Yu.
\newblock Convergence and superconvergence of staggered discontinuous
  {G}alerkin methods for the three-dimensional maxwell's equations on cartesian
  grids.
\newblock {\em J. Comput. Phys.}, 235:14--31, 2013.

\bibitem{OLS}
E.~Chung, H.~Kim, and O.~Widlund.
\newblock Two-level overlapping schwarz algorithms for a staggered
  discontinuous {G}alerkin method.
\newblock {\em SIAM J. Numer. Anal.}, 51:47--67, 2013.

\bibitem{SDG-cd}
E.~Chung and C.~Lee.
\newblock A staggered discontinuous {G}alerkin method for the
  convection-diffusion equation.
\newblock {\em J. Numer. Math.}, 20:1--31, 2012.

\bibitem{SDG-curl}
E.~Chung and C.~Lee.
\newblock A staggered discontinuous {G}alerkin method for the curl-curl
  operator.
\newblock {\em IMA J. Numer. Anal.}, 32:1241--1265, 2012.

\bibitem{CiCP}
E.~Chung and W.~Leung.
\newblock A sub-grid structure enhanced discontinuous {G}alerkin method for
  multiscale diffusion and convection-diffusion problems.
\newblock {\em Commum. Comput. Phys.}, 14:370--392, 2013.

\bibitem{basabe:562}
Jonas~D. {De Basabe} and Mrinal~K. Sen.
\newblock New developments in the finite-element method for seismic modeling.
\newblock {\em The Leading Edge}, 28(5):562--567, 2009.

\bibitem{delprat-jannaud:T37}
Florence Delprat-Jannaud and Patrick Lailly.
\newblock Wave propagation in heterogeneous media: Effects of fine-scale
  heterogeneity.
\newblock {\em Geophysics}, 73(3):T37--T49, 2008.

\bibitem{egh12}
Y.~Efendiev, J.~Galvis, and T.~Hou.
\newblock Generalized multiscale finite element method.
\newblock 2012.

\bibitem{eglms13}
Y.~Efendiev, J.~Galvis, R.~Lazarov, M.~Moon, and M.~Sarkis.
\newblock Generalized multiscale finite element method. \uppercase{IPDG}.
\newblock Submitted, 2013.

\bibitem{eglp13}
Y.~Efendiev, J.~Galvis, G.~Li, and M.~Presho.
\newblock Generalized multiscale finite element methods. oversampling
  strategies.
\newblock nternational Journal for Multiscale Computational Engineering, to
  appear, 2013.

\bibitem{Geophysics}
R.~Gibson, K.~Gao, E.~Chung, and Y.~Efendiev.
\newblock Multiscale modeling of acoustic wave propagation in two-dimensional
  media.
\newblock {\em Under Review}.

\bibitem{IPDG}
M.~Grote and D.~Schotzau.
\newblock Optimal error estimates for the fully discrete interior penalty {DG}
  method for the wave equation.
\newblock {\em J. Sci. Comput.}, 40:257--272, 2009.

\bibitem{Hermann:2011fk}
Verena Hermann, Martin K{\"a}ser, and Crist{\'o}bal~E. Castro.
\newblock Non-conforming hybrid meshes for efficient {2-D} wave propagation
  using the discontinuous galerkin method.
\newblock {\em Geophysical Journal International}, 184(2):746--758, 2011.

\bibitem{kaser:76}
Martin K\"{a}ser, Christian Pelties, Cristobal~E. Castro, Hugues Djikpesse, and
  Michael Prange.
\newblock Wavefield modeling in exploration seismology using the discontinuous
  {Galerkin} finite-element method on hpc infrastructure.
\newblock {\em The Leading Edge}, 29(1):76--85, 2010.

\bibitem{springerlink:10.1007/s00450-010-0109-1}
Dimitri Komatitsch, Dominik G{\"o}ddeke, Gordon Erlebacher, and David
  Mich{\'e}a.
\newblock Modeling the propagation of elastic waves using spectral elements on
  a cluster of 192 gpus.
\newblock {\em Computer Science - Research and Development}, 25:75--82, 2010.
\newblock 10.1007/s00450-010-0109-1.

\bibitem{GJI:GJI967}
Dimitri Komatitsch and Jeroen Tromp.
\newblock Introduction to the spectral element method for three-dimensional
  seismic wave propagation.
\newblock {\em Geophysical Journal International}, 139(3):806--822, 1999.

\bibitem{GJI:GJI1653}
Dimitri Komatitsch and Jeroen Tromp.
\newblock Spectral-element simulations of global seismic wave propagation---i.
  validation.
\newblock {\em Geophysical Journal International}, 149(2):390--412, 2002.

\bibitem{GJI:GJI3620}
B.~Lombard, J.~Piraux, C.~{G\'elis}, and J.~Virieux.
\newblock Free and smooth boundaries in {2-D} finite-difference schemes for
  transient elastic waves.
\newblock {\em Geophysical Journal International}, 172(1):252--261, 2008.

\bibitem{masson:N33}
Y.~J. Masson and S.~R. Pride.
\newblock Finite-difference modeling of {Biot's} poroelastic equations across
  all frequencies.
\newblock {\em Geophysics}, 75(2):N33--N41, 2010.

\bibitem{GJI:GJI5221}
Peter Moczo, Jozef Kristek, Martin Galis, Emmanuel Chaljub, and Vincent
  Etienne.
\newblock {3-D} finite-difference, finite-element, discontinuous-{Galerkin} and
  spectral-element schemes analysed for their accuracy with respect to p-wave
  to s-wave speed ratio.
\newblock {\em Geophysical Journal International}, 187(3):1645--1667, 2011.

\bibitem{GJI:GJI4985}
Christina Morency, Yang Luo, and Jeroen Tromp.
\newblock Acoustic, elastic and poroelastic simulations of {CO$_2$}
  sequestration crosswell monitoring based on spectral-element and adjoint
  methods.
\newblock {\em Geophysical Journal International}, pages no--no, 2011.

\bibitem{Pelties:2010kx}
C.~Pelties, M.~K{\"a}ser, V.~Hermann, and C.~E. Castro.
\newblock Regular versus irregular meshing for complicated models and their
  effect on synthetic seismograms.
\newblock {\em Geophysical Journal International}, pages no--no, 2010.

\bibitem{IPDGbook}
Beatrice~M. Riviere.
\newblock {\em Discontinuous {G}alerkin Methods For Solving Elliptic And
  parabolic Equations: Theory and Implementation}.
\newblock SIAM, 2008.

\bibitem{saenger:SM293}
Erik~H. Saenger, Radim Ciz, Oliver~S. Kr\"{u}ger, Stefan~M. Schmalholz, Boris
  Gurevich, and Serge~A. Shapiro.
\newblock Finite-difference modeling of wave propagation on microscale: A
  snapshot of the work in progress.
\newblock {\em Geophysics}, 72(5):SM293--SM300, 2007.

\bibitem{symes:2602}
William Symes, Igor~S. Terentyev, and Tetyana Vdovina.
\newblock Getting it right without knowing the answer: Quality control in a
  large seismic modeling project.
\newblock {\em SEG Technical Program Expanded Abstracts}, 28(1):2602--2606,
  2009.

\bibitem{virieux:1933}
Jean Virieux.
\newblock {SH}-wave propagation in heterogeneous media: Velocity-stress
  finite-difference method.
\newblock {\em Geophysics}, 49(11):1933--1942, 1984.

\bibitem{virieux:889}
Jean Virieux.
\newblock {P-SV} wave propagation in heterogeneous media: Velocity-stress
  finite-difference method.
\newblock {\em Geophysics}, 51(4):889--901, 1986.

\bibitem{CMS}
W.~Zhang, L.~Tong, and E.~Chung.
\newblock Exact nonreflecting boundary conditions for three dimensional
  poroelastic wave equations.
\newblock {\em Comm. Math. Sci.}, To appear.

\bibitem{JCAM-lod}
W.~Zhang, L.~Tong, and E.~Chung.
\newblock A new high accuracy locally one-dimensional scheme for the wave
  equation.
\newblock {\em J. Comput. Appl. Math.}, 236:1343--1353, 2011.

\bibitem{nmtma}
W.~Zhang, L.~Tong, and E.~Chung.
\newblock Efficient simulation of wave propagation with implicit finite
  difference schemes.
\newblock {\em Numer. Math. Theor. Meth. Appl.}, 5:205--228, 2012.

\end{thebibliography}

\end{document}